\documentclass[11pt]{amsart}
\usepackage[french,english]{babel}
\usepackage{amsmath}
\usepackage{amsfonts}
\usepackage{amssymb}
\usepackage{amsthm}
\usepackage{graphicx}
\usepackage[latin1]{inputenc}
\usepackage[T1]{fontenc}
\usepackage{eufrak}
\usepackage[all]{xy}
\newtheorem{thm}{Theorem}[section]
\theoremstyle{definition}

\theoremstyle{remark}
\newtheorem{rem}{Remark}
\usepackage{color}
\newenvironment{dem}{\noindent{\textit{Proof.}}}{\begin{flushright}$\diamondsuit$\end{flushright}}
\theoremstyle{plain}
\newtheorem{lem}{\bf{Lemma}}[section]
\theoremstyle{plain}
\newtheorem{prop}{Proposition}[section]
\newtheorem{cor}{Corollary}
\numberwithin{equation}{subsection}
\begin{document}
\author{Guillaume Pouchin}
\def\sym{{\mathfrak{S}}}
\def\Zz{{\mathbb{Z}}}
\def\Nn{{\mathbb{N}}}
\def\Cc{{\mathbb{C}}}
\def\Qq{{\mathbb{Q}}}
\def\End{{\mathrm{End}}}
\def\Hom{{\mathrm{Hom}}}
\def\Im{{\mathrm{Im}\,}}
\def\Aut{{\mathrm{Aut}}}
\def\Ker{{\mathrm{Ker}\,}}
\def\dim{{\mathrm{dim}\,}}
\def\Ext{{\mathrm{Ext}^1}}
\def\P1{{\mathbb{P}^1}}
\def\CalO{{\mathcal{O}}}
\def\CalF{{\mathcal{F}}}
\def\CalG{{\mathcal{G}}}
\def\CalE{{\mathcal{E}}}
\def\CalH{{\mathcal{H}}}
\def\CalI{{\mathcal{I}}}
\def\Hilb{{\mathrm{Hilb}}}
\def\Coh{{\underline{\mathrm{Coh}}}}
\def\Vect{{\mathrm{Vect}}}
\def\Vec{{\underline{\mathrm{Bun}}}}
\def\Tor{{\underline{\mathrm{Tor}}}}
\def\Lie{{\mathfrak{g}}}
\def\Higgsa{{\underline{\Lambda}^\alpha_{\mathbb{P}^1}}}
\def\Higgsb{{\underline{\Lambda}^\beta_{\mathbb{P}^1}}}
\def\Higgsab{{\underline{\Lambda}^{\alpha+\beta}_{\mathbb{P}^1}}}
\def\Higgs{{\underline{\Lambda}_{\mathbb{P}^1}}}
\def\alg{{\mathcal{H}}}
\def\alga{{\mathcal{H}^\alpha}}
\title{Higgs algebra of curves and loop crystals}

\maketitle
\begin{abstract}
We define the Higgs algebra $\mathcal{H}_\P1$ of the projective line,
as a convolution algebra of constructible functions on the global
nilpotent cone $\underline{\Lambda}_\P1$, a lagrangian substack of the
Higgs bundle $T^*\Coh_\P1$, where $\Coh_\P1$ is the
  stack of coherent sheaves on $\P1$. We prove that $\mathcal{H}_\P1$
  is isomorphic to (some completion of) $U^+(\widehat{sl}_2)$. We use
  this geometric realization to define a semicanonical basis of
  $U^+(\widehat{sl}_2)$, indexed by irreducible components of
  $\underline{\Lambda}_\P1$. We also construct a combinatorial data on
  this set of irreducible components in the spirit of \cite{KS}, which
  is an affine analog of a crystal. We call it a loop crystal and give
  some of its properties.
\end{abstract}
\section*{introduction}
The link between Kac-Moody algebra and geometry of quiver
representation is known since the work of Ringel (\cite{Ri}) and Lusztig (\cite{L2}), and has lead to the theory of canonical bases of quantum algebras. More recently, Kapranov and Bauman-Kassel (\cite{BK}) have considered the Hall algebra of the category of coherent sheaves on the projective line and showed that some composition subalgebra is isomorphic to a positive part of the affine quantum algebra of $sl_2$. Schiffmann (\cite{Sc2}) generalizes this situation to any weighted projective line, and showed that the algebras obtained are affine versions of some quantum algebras. He used this setting to define canonical bases for some of these algebras.

In the context of quivers, Lusztig considered also a geometric
construction of (enveloping) Kac-moody algebras via constructible
functions on some nilpotent part $\Lambda$ of the cotangent bundle of
the stack of representations of a quiver. He then defined
semicanonical bases using that geometry, whose elements are indexed by irreducible components of $\Lambda$. Kashiwara and Saito (\cite{KS}) later considered natural correspondences between irreducible components to provide a geometric construction of the crystal associated to the semicanonical bases. 

This article is an attempt to develop Kashiwara and Saito's ideas in
the context of curves. We define the Higgs algebra of a smooth
projective curve $X$ as constructible functions on the Hitchin stack $\underline{\Lambda}$, which is a lagrangian substack of the Higgs bundle (i.e. the cotangent bundle of the stack of coherent sheaves on $X$). We prove that if $X= \P1$ the Higgs algebra is isomorphic to (some completion of) the positive part (in Drinfeld sense) of the enveloping algebra of $\widehat{sl}_2$, and define a semicanonical basis indexed by irreducible components Irr($\Higgs$). We then construct a stratification of the stack $\Higgs$ and define nice correspondences between strata which give rise  to natural operators acting on the set Irr($\Higgs$). We call loop crystal this new combinatorial data.

We then describe some properties of our loop crystal. First it generalizes the notion of (affine) crystal as it contains a crystal data, but also many more operators. Unlike the $\widehat{sl}_2$-crystal in it, this loop crystal is connected. Moreover, if we follow ideas from quiver varieties and consider some stable conditions on the Hitchin stack, the subcrystal looks like some limit of Kirillov-Reshetikhin modules, which are (conjecturally) the only finite dimensional modules with crystal bases.
This suggests that the notion of loop crystal could be axiomatized from the case of $\widehat{sl}_2$ (as the notion of crystal was from $sl_2$) and that representation theoric information may be recovered from it.

We should also mention that our construction are still valid in the
case of weighted projective lines. The operators obtained in the
combinatorial data are then indexed by simple rigid objects. However,
as the positive part of the enveloping algebra is not the Drinfeld
part, the link with representation theory is less clear. We plan to
develop this in a further paper.

The paper is organized as follows: the first part is about the geometry of the stack of coherent sheaves on $\P1$ and the definition of the stack $\underline{\Lambda}$, which is the analog of the nilpotent cone. In the second part we describe the irreducible component of $\underline{\Lambda}$ in the case of $\P1$. The third part contain the definition of the Higgs algebra $\alg$ and its first properties. Then we in the next section we prove the existence (and unicity) of a semicanonical basis of $\alg$. We identify the algebra $\alg$ with the positive part of the enveloping algebra of $\widehat{sl_2}$ in the fifth part. The last section is devoted to the loop crystal structure on the set of irreducible components of $\underline{\Lambda}$.

\textbf{Acknowledgements}. 
This article is part of the author's Phd thesis at Paris VI University
under the supervision of O. Schiffmann. He deeply thanks him for his
constant support, encouragement and availability.
\section{The global analog of the nilpotent cone}
\subsection{}
Let $\text{Coh}_\P1$ be the category of coherent sheaves on $\P1(\Cc)$.\\
We write $K(\text{Coh}_\P1)$ for the Grothendieck group and $K^+(\text{Coh}_\P1)$ for the
positive cone, which consists of elements $\alpha$ for which there exists some
sheaf $\CalF$ with $[\CalF]=\alpha$. It is isomorphic to $(\Nn^*
\times \Zz) \cup (\{0\}\times \Nn^*)$, and the
isomorphism is given by the rank and the degree of a sheaf. The
Grothendieck group is equipped with the Euler form $\langle [M],[N]\rangle=\dim
\Hom(M,N)- \dim
\Ext(M,N)=\text{rk}(M)\text{rk}(N)+\text{rk}(M)\text{deg}(N)-\text{deg}(M)\text{rk}(N)$,
where $\text{rk}(M)$ and $\text{deg}(N)$ are the rank and the degree
of $N$.\\
For $\alpha=(r,d) \in K^+(\text{Coh}_\P1)$, let $\Coh_\P1^\alpha$ denote the stack
classifying coherent sheaves of class $\alpha$ on $\P1$. It is known
(see \cite{LaMB}) that $\Coh_\P1^\alpha$ is smooth and connected of
dimension $-\langle \alpha,\alpha\rangle=-r^2$.\\
 This stack has a local presentation as follows (see
 \cite{Gr},\cite{Le} or \cite{Sc1}).\\
Let $\alpha \in K^+(\text{Coh}_\P1)$ and $\CalE \in \Coh_\P1$.
Define the following functor from the category of affine schemes over $\Cc$ to the
category of sets:
\begin{align*}
\text{\underline{Hilb}}_{\CalE,\alpha}(\Sigma)=\{ \phi_\Sigma: \CalE
\boxtimes \CalO_\Sigma \twoheadrightarrow
\CalF, \ \CalF \text{ is a coherent $\Sigma$-flat sheaf, } & \\
& \hspace{-6cm} \CalF_\sigma\text{ is of class $\alpha$ for all closed point $\sigma \in \Sigma$} \}/\sim,
\end{align*}
where two such morphisms are equivalent if they have the same kernel.\\
This functor is representable by a projective scheme
$\Hilb_{\CalE,\alpha}$ (see \cite{Gr}).\\
Fix an integer $n\in \Nn$ and define $d(n,
\alpha)=\langle[\CalO(n)],\alpha\rangle$ and
$\CalE_n^\alpha=\Cc^{d(n,\alpha)}\otimes \CalO(n)$ if $d(n,\alpha)\geq
0$. A map $\phi_\Sigma$ induces
for each closed point $\sigma \in \Sigma$ a linear map
$\phi_{*,\sigma}:\Cc^{d(n,\alpha)} \rightarrow
\Hom(\CalO(n),\CalF_\sigma$).\\
Let us consider the subfunctor defined by:
\begin{align*}
\Sigma \mapsto \{ (\phi_\Sigma: \CalE_n^\alpha \boxtimes \CalO_\Sigma \twoheadrightarrow
\CalF ) \in \text{\underline{Hilb}}_{\CalE_n^\alpha,\alpha}(\Sigma), \ | \forall
\sigma \in \Sigma, & \\ &
\hspace{-3cm} \phi_{*,\sigma} :\Cc^{d(n,\alpha)} \simeq
\Hom(\CalO(n),\CalF_\sigma) \}/\sim.
\end{align*}
This subfunctor is representable by a smooth open quasiprojective
subscheme $Q_n^\alpha$ of $\Hilb_{\CalE_n^\alpha,\alpha}$ (see \cite{Le}).\\
The group $G_n^\alpha=\Aut(\CalE_n^\alpha)=GL(d(n,\alpha))$ acts
naturally on $\Hilb_{\CalE_n^\alpha,\alpha}$ and $Q_n^\alpha$.\\
The quotient stacks $\Coh_\P1^{\alpha,\geq n}:=[Q_n^\alpha/G_n^\alpha]$ for $n \in \Zz$ are open
substacks of $\Coh_\P1^\alpha$ which form
an atlas, i.e.: 
\[
\Coh_\P1^\alpha=\bigcup_{n \in \Zz}\Coh_\P1^{\alpha,\geq
  n}=\bigcup_{n\in \Zz} [Q_n^\alpha/G_n^\alpha]
\]
We also introduce the stack $\Vec_\P1=\cup_\alpha \Vec_\P1^\alpha$ of locally free sheaves on
$\P1$, which is an open substack of $\Coh_\P1$. We have an atlas given by the open substacks
$\Vec_\P1^{\alpha,\geq n}=[U_n^\alpha/G_n^\alpha]$, where
$U_n^\alpha=\{ (\phi:\CalE_n^\alpha \twoheadrightarrow \CalF)\in Q_n^\alpha, \CalF \text{ is locally free}\}/\sim$.
\subsection{}
We want to describe the cotangent stack $T^*(\Coh_\P1^\alpha)$ by giving it an
atlas obtained by symplectic reduction of the varieties
$T^*(Q_n^\alpha)$.\\
First recall that the tangent space at a point $\phi:\CalE
\twoheadrightarrow \CalF$ of $Q_n^\alpha$ is canonically isomorphic to
$\Hom(\Ker \phi,\CalF)$ (see \cite{Le}).\\
The group $G_n^\alpha$ acts on $T^*Q_n^\alpha$ in a Hamiltonian fashion.
The corresponding moment map $\mu_n: T^*Q_n^\alpha \rightarrow
(\Lie_n^\alpha)^*$ is described as follows: 
 over a point $z=(\phi,f)\in T^*Q_n^\alpha$ with $\phi:
 \CalO(n)^{ d(n,\alpha)} \twoheadrightarrow \CalF$ and $f\in \Hom(\Ker
 \phi,\CalF)^*$, we have 
\[
\begin{array}{cccc}
\mu_n^\alpha(z): &\Lie_n^\alpha&\rightarrow &\Cc\\
 & g  & \mapsto & <f,(\phi \circ g)|_{\Ker \phi}>,
\end{array}
\]
where $\Lie$ acts on $\Hom(\CalO(n)^{d(n,\alpha)},\CalF)$ by means
of the isomorphism $\phi_*: \Cc^{d(n,\alpha)} \simeq
\Hom(\CalO(n),\CalF)$.\\
We want to describe the subvariety $(\mu_n^\alpha)^{-1}(0) \subseteq
T^*Q_n^\alpha$. To do this, fix some point $\phi: \CalO(n)^{d(n,\alpha)}
\twoheadrightarrow \CalF$ in $Q_n^\alpha$ and write the short exact sequence:
\[
0 \rightarrow \Ker(\phi) \rightarrow \CalO(n)^{\oplus d(n,\alpha)}
\rightarrow \CalF \rightarrow 0.
\]
Apply the functor $\Hom(\_,\CalF)$:
\begin{multline*}
0 \rightarrow \Hom(\CalF,\CalF)\rightarrow
\Hom(\CalO(n)^{ d(n,\alpha)},\CalF) \rightarrow \Hom(\Ker\phi, \CalF) \\
\rightarrow \Ext(\CalF,\CalF)  \rightarrow
\Ext(\CalO(n)^{ d(n,\alpha)},\CalF) \rightarrow \cdots
\end{multline*}
Since $\phi$ belongs to $Q_n^\alpha$, $\langle\CalO(n),\CalF \rangle=d(n,\alpha)=\dim\Hom(\CalO(n),\CalF)$ so we have
$\Ext(\CalO(n),\CalF)=0$. Dualizing, we get :
\[
0 \rightarrow \Ext(\CalF,\CalF)^* \rightarrow \Hom(\Ker\phi,\CalF)^*
\xrightarrow{a} \Hom(\CalO(n)^{d(n,\alpha)},\CalF)^* \rightarrow \cdots
\]
One checks that the map $a$ is precisely the moment map $\mu_n$. So if
$(\phi,f)$ is in $(\mu_n^\alpha)^{-1}(0)$ then $f$ defines a unique
element in $\Ext(\CalF,\CalF)^*$, which we still denote $f$.\\
Serre duality gives a canonical isomorphism: 
$\Ext(\CalF,\CalF)^* \simeq \Hom(\CalF,\CalF(-2))$, where we write
$\CalF(-2)$ for $\CalF \otimes \CalO(-2)$.\\
We finally have :
\[
(\mu_n^\alpha)^{-1}(0)^\alpha=\{ (\phi : \CalO(n)^{d(n,\alpha)} \twoheadrightarrow
\CalF,f) \in T^*Q_n^\alpha\  \vert \ f \in \Hom(\CalF,\CalF(-2))\}/\sim
\]
By symplectic reduction the cotangent bundle stack of the
quotient stack $Q_n^\alpha/G_n^\alpha$ is the quotient
$[(\mu_n^\alpha)^{-1}(0)/G_n^\alpha]$. This gives us an atlas of
$T^*(\Coh_\P1)$:
\[
T^*(\Coh_\P1^\alpha)=\bigcup_{n\in \Zz}
[(\mu_n^\alpha)^{-1}(0)/G_n^\alpha]
\]
The cotangent stack $T^*(\Coh_\P1^\alpha)$ represents the functor
$\underline{\text{Higgs}}_\P1^\alpha$ from the category of affine
schemes over $\Cc$ to the category of groupoids, where we write
$\CalO_{\Sigma\times\P1}(-2)$ for  $\CalO_\Sigma \boxtimes \CalO_\P1(-2)$:
\begin{align*}
\underline{\text{Higgs}}_\P1^\alpha(\Sigma)=\{ (\CalF,f), \ \CalF \text{ is a
  coherent $\Sigma$-flat sheaf on $\P1\times \Sigma$, } & \\
& \hspace{-6cm}\CalF_\sigma\text{ is of class $\alpha$ for all closed
  point $\sigma \in \Sigma$}, \\
& \hspace{-4.5cm}f \in \Hom(\CalF, \CalF \otimes \CalO_{\Sigma \times \P1}(-2))\}
\end{align*}
where a morphism $\psi$ between two objects $(\CalF,f)$ and
$(\CalF',f')$ is an isomorphism $\psi: \CalF \simeq \CalF'$ such that
the following diagram commute:
\[
\xymatrix{
\CalF \ar[r]^\psi \ar[d]^f & \CalF' \ar[d]^{f'}\\
\CalF\otimes\CalO_{\Sigma \times \P1}(-2) \ar[r]^{\psi} & \CalF'\otimes\CalO_{\Sigma \times \P1}(-2)
}
\]
\subsection{}
Let us now introduce the nilpotent part:
\[
S_n^\alpha:=\mu_n^{-1}(0)^{\alpha,\text{nilp}}=\{
(\phi:\CalO(n)^{d(n,\alpha)} \twoheadrightarrow \CalF,f)\in \mu_n^{-1}(0), \ f \text{
nilpotent}\}
\]
where we say that $f$ is nilpotent if there exists $m$ such that
\[
f(-2(m-1))\circ\cdots f(-2)\circ f =0
\]
as an element of $\Hom(\CalF,\CalF(-2m))$.\\
The quotient stacks $\underline{\Lambda}_\P1^{\alpha,\geq
  n}=[S_n^\alpha/G_n^\alpha]$ are closed substacks of
  $T^*\Coh_\P1^{\alpha,\geq n}$, and form a compatible family with
  respect to the inductive system
  $T^*\Coh_\P1^{\alpha,\geq n}$. They give rise in the
  limit to a closed substack
\[\underline{\Lambda}_\P1^\alpha=\varinjlim
  [S_n^\alpha/G_n^\alpha]=\bigcup_{n\in \Zz} [S_n^\alpha/G_n^\alpha]
  \subseteq T^*\Coh_\P1^\alpha
\]
The stack $\Higgsa$ represents the functor
$\underline{\text{Higgs}}_\P1^{\alpha,\text{nilp}}$ from the category of affine schemes
over $\Cc$ to the category of groupoids:
\begin{align*}
\underline{\text{Higgs}}_\P1^{\alpha,\text{nilp}}(\Sigma)= \{ (\CalF,f), \ \CalF \text{ is a
  coherent $\Sigma$-flat sheaf on $\P1 \times \Sigma$, } & \\
& \hspace{-6cm}\CalF_\sigma\text{ is of class $\alpha$ for all closed
  point $\sigma \in \Sigma$}, \\
 & \hspace{-5.5cm}f \in \Hom(\CalF, \CalF \otimes \CalO_{\Sigma\times\P1}(-2))\text{ is nilpotent}\}
\end{align*}
where the morphism between two objects are the same as for the functor $\underline{\text{Higgs}}_\P1^\alpha$.\\
We also have the same description for the stack $T^*\Vec_\P1$. In that
case the nilpotency condition for $f$ is empty. Indeed, for every
vector bundle $\mathcal{V}$ we have
$\Hom(\mathcal{V},\mathcal{V}(-2k))=0$ for $k>>0$. This implies that
$T^*\Vec_\P1$ is also an open substack of $\Higgs$. We write
$T^*\Vec_\P1^{\alpha,\geq n}= [R_n^\alpha/G_n^\alpha]$ where
$R_n^\alpha=\{(\phi:\CalE^\alpha_n\twoheadrightarrow\CalF,f)\in S_n^\alpha,\CalF \text{ is locally free}\}$.
\section{Irreducible components of $\Higgsa$}\label{Irredcomp}
In this section we want to describe the irreducible components of
$\Higgsa$.
\subsection{}
We begin with a lemma for the irreducible components of the torsion
part. For a partition $\lambda$ of $d$ denote by $\CalO_\lambda$ the
smooth strata of $\Coh_\P1^{(0,d)}$ parametrizing the sheaves $\{
\CalO_{x_1}^{(\lambda_1)}\oplus \cdots \oplus
\CalO_{x_d}^{(\lambda_d)}| x_i \neq x_j\}$, where $\CalO_x^{(d)}$ is the
indecomposable torsion sheaf of degree $d$ supported at $x$. Let
$T_{\CalO_\lambda}^*\Coh_\P1^{(0,d)}$ be the conormal bundle to this strata.
\begin{lem}[\cite{La1},theorem 3.3.13]\label{lemmetorsion}
We have the decomposition into irreducible components:
\[
\underline{\Lambda}_\P1^{(0,d)}=\bigcup_{\lambda \vdash d}
\overline{T_{\CalO_\lambda}^*\Coh_\P1^{(0,d)}}
\]
and each has (stacky) dimension $0$.
\end{lem}
Now for a locally free sheaf $\mathcal{V}$ of rank $r$ and
degree $d'$ and a partition $\lambda$ of $d''$ with $d=d'+d''$ let
$X_{\mathcal{V},\lambda}$ be the substack
 of $\Higgsa$ parametrizing pairs $(
\mathcal{V} \oplus \tau,f)$ with $\tau \in \CalO_\lambda$.
\begin{thm}\label{irr}
The irreducible components of $\underline{\Lambda}_\P1^{(r,d)}$ are exactly:
\[
\text{Irr}(\underline{\Lambda}_\P1^{(r,d)})=\{ \overline{X}_{\mathcal{V},\lambda} \}_{\mathcal{V},\lambda}
\]
Each is of (stacky) dimension $-r^2$.
\end{thm}
To prove this theorem, we will proceed in two steps. In the first
step, we stratify our space into locally closed subspaces where the degree of
the torsion part of $\CalF$ is fixed. In the second step, we split
every strata between its locally free and torsion part. Denote
$\CalF^\text{tor}$ the torsion part of $\CalF$. It is a subsheaf of $\CalF$.\\
For $l \in \Nn$, define a locally closed stack of $\Coh_\P1^{(r,d)}$
which parametrizes isomorphism classes of objects:
\[
\Coh_\P1^{r,d,l}=\{ \CalF \in \Coh_\P1^{r,d}\quad | \quad \text{deg}(\CalF^\text{tor})=l\}
\]
We have
\[
\Coh_\P1^{r,d}=\bigsqcup_{l \in \Nn}\Coh_\P1^{r,d,l} \text{  (a
  locally finite union)}
\]
Denoting by $\pi: \Higgsa \rightarrow \Coh_\P1^\alpha$ the first
projection, we set $\underline{\Lambda}_\P1^{\alpha,l}=\pi^{-1}(\Coh_\P1^{\alpha,l})$.\\
For each irreducible component $Z$ of $\Higgsa$ there is a unique
integer $l$ such that $Z\cap \underline{\Lambda}_\P1^{r,d,l}$ is dense
in $Z$. We start by describing the
irreducible components of $\underline{\Lambda}_\P1^{r,d,l}$.

\subsection{}
Denote by $\underline{L}^{r,d,l}$ the stack parametrizing isomorphism
classes of objects:
\begin{align*}
\underline{L}^{r,d,l}:=\{(\mathcal{V},\tau,f_1,f_2,f_3),\mathcal{V} \in
\Vec_\P1^{r,d-l},\tau \in \Coh_\P1^{(0,d)},f_1 \in
\Hom(\mathcal{V},\mathcal{V}(-2)),\\f_2\in\Hom(\tau,\tau) \text{
  nilpotent},f_3\in\Hom(\mathcal{V},\tau)\}
\end{align*}
where a morphism $\psi$ between two objects $(\mathcal{V},\tau,
f_1,f_2,f_3)$ and $(\mathcal{V}',\tau',f_1',f_2',f_3')$ is just a
couple $(\psi_1,\psi_2)$ with $\psi_1: \mathcal{V} \simeq
\mathcal{V}'$ and $\psi_2: \tau \simeq \tau'$ such that the three
following diagrams commute:
\[
\xymatrix{
\mathcal{V} \ar[r]^{\psi_1} \ar[d]^{f_1} &\mathcal{V}' \ar[d]^{f_1'} &
\tau \ar[r]^{\psi_2} \ar[d]^{f_2}& \tau' \ar[d]^{f_2'} & \mathcal{V}
\ar[r]^{\psi_1} \ar[d]^{f_3} & \mathcal{V}'\ar[d]^{f_3'}\\
\mathcal{V} \ar[r]^{\psi_1} & \mathcal{V}'& \tau \ar[r]^{\psi_2} &
\tau' & \tau \ar[r]^{\psi_2} & \tau'
}
\]
We have a natural diagram:
\begin{equation}\label{Dec}
\xymatrix{
 & \underline{L}^{r,d,l} \ar[dl]_{\pi_1} \ar[dr]^{\pi_2} & \\
\underline{\Lambda}_\P1^{r,d,l} & & T^*\Vec^{r,d-l} \times \underline{\Lambda}_\P1^{0,l}
}
\end{equation}
where $\pi_1$ is defined from the functor of groupoids:
\[
\begin{array}{cccc}
\pi_1: &(\mathcal{V},\tau,f_1,f_2,f_3)& \mapsto & \Bigg( \mathcal{V}
\oplus \tau,\begin{pmatrix}f_1&0\\f_3&f_2\end{pmatrix}\Bigg) \\
 & (\psi_1,\psi_2) & \mapsto & \begin{pmatrix}\psi_1 & 0 \\ 0 & \psi_2 \end{pmatrix}
\end{array}
\]
and for $\pi_2$:
\[
\begin{array}{cccc}
\pi_2: & (\mathcal{V},\tau,f_1,f_2,f_3) & \mapsto &
((\mathcal{V},f_1),(\tau,f_2)) \\
 & (\psi_1,\psi_2)& \mapsto & (\psi_1,\psi_2)
\end{array}
\]
\begin{lem}\label{buntor}
The map $\pi_1$ is an affine fibration and $\pi_2$ is a vector bundle,
each is of relative dimension $lr$ and with connected fibers. This
induces a bijection between irreducible components:
\[
Irr(\underline{\Lambda}_\P1^{r,d,l}) \leftrightarrow
Irr(T^*\Vec_\P1^{r,d-l}) \times
Irr(\underline{\Lambda}_\P1^{0,l}).
\]
Moreover this correspondence preserves dimensions, i.e. if we have $Z
\leftrightarrow Z_1 \times Z_2$ under this correspondence, then $\dim
Z=\dim Z_1 + \dim Z_2$.
\end{lem}
\begin{dem}
First we recall that $\pi_1$ is well defined because the
nilpotency condition is empty for $f_1$, so that $f$ is indeed nilpotent.\\
The result is obvious for $\pi_2$ since $\dim\Hom(\mathcal{V},\tau)=rl$.\\
We introduce the following natural stack:
\[
\underline{\mathcal{S}}^{r,d,l}= \{ (\CalF,f)| \CalF \in \Coh_\P1^{r,d,l},\ f\in
\End(\CalF), f|_{\CalF^\text{tor}}=0, f(\CalF)\subseteq
\CalF^\text{tor} \}
\]
with a morphism $\psi$ between objects $(\CalF,f)$ and $(\CalF',f')$
is an isomorphism $\psi:\CalF \simeq \CalF'$ such that the diagram
\[
\xymatrix{
 \CalF \ar[r]^\psi \ar[d]^f & \CalF' \ar[d]^{f'} \\
\CalF \ar[r]^\psi & \CalF'
}
\]
is commutative.\\
We have a natural map $\pi: \underline{\mathcal{S}}^{r,d,l}
\rightarrow \Coh_\P1^{r,d,l}$, which makes
$\underline{\mathcal{S}}^{r,d,l}$ a vector bundle over $\Coh_\P1^{r,d,l}$.\\
Define its pullback over $\underline{\Lambda}_\P1^{r,d,l}$:
\[
\underline{\tilde{\mathcal{S}}}^{r,d,l}=\underline{\mathcal{S}}^{r,d,l}\times_{\Coh_\P1^{r,d,l}}\underline{\Lambda}_\P1^{r,d,l}
\]
This is a vector bundle over $\underline{\Lambda}_\P1^{r,d,l}$ of rank
$rl$.\\
We have a natural action of $\underline{\tilde{\mathcal{S}}}^{r,d,l}$
on $\underline{L}^{r,d,l}$ over $\underline{\Lambda}_\P1^{r,d,l}$
defined as follows. Take a point $P=(\mathcal{V},\tau, f_1,f_2,f_3) \in
\underline{\Lambda}_\P1^{r,d,l}$. The fiber of
$\underline{\tilde{\mathcal{S}}}^{r,d,l}$ over $\pi_1(P)$ is by construction
canonically identified with $\{ g \in \End(\mathcal{V} \oplus \tau)|
g(\tau)=0,\ g(\mathcal{V} \oplus \tau) \subseteq \tau \}$. We define
the action as follows:
\[
g.P=(\mathcal{V},\tau,f_1,f_2,f_3-gf_1+f_2g)
\]
which corresponds to the action of $(Id +g)$ by conjugation on
$\begin{pmatrix} f_1 & 0 \\ f_3 & f_2 \end{pmatrix}$.\\
As we have $\Aut(\mathcal{V}\oplus \tau)= (\Aut(\mathcal{V}) \times
\Aut(\tau))\rtimes \Hom(\mathcal{V},\tau)$, we can identify
$\underline{\Lambda}_\P1^{r,d,l}$ as the quotient of
$\underline{L}^{r,d,l}$ by the action of
$\underline{\tilde{S}}^{r,d,l}$.
\end{dem}
It remains to describe the irreducible components of
$T^*\Vec_\P1^{(r,k)}$.
\begin{lem}\label{IrrBun}
The irreducible components of $T^*\Vec^{(r,k)}$ are the closures of the
conormal bundles $T^*_\mathcal{V}\Vec^{(r,k)}$, for $\mathcal{V} \in
\Vec_\P1^{(r,k)}$. Each is of dimension $\dim \Vec_\P1^{(r,k)}=-r^2$.
\end{lem}
\begin{dem}
We have $T^*\Vec^{(r,k)}=\varinjlim T^*\Vec_\P1^{(r,k),\geq
  n}=\varinjlim [R_n^{(r,k)}/G_n^{(r,k)}]$. Since $R_n^{(r,k)}$ has
  finitely many $G_n^{(r,k)}$ orbits $\CalO$ we deduce that
  $R_n^{(r,k)}=T^*\Vec^{(r,k),\geq n} \cap \bigcup_{\CalO \subseteq U_n^{(r,k)}}
  T_\CalO^*R_n^{(r,k)}$ has $\{\overline{T_\CalO^*U_n^{(r,k)}}\}$ as
  irreducible components, each of the same dimension. Hence $T^*\Vec_\P1^{(r,k),\geq
  n}$ has $\{\overline{T_\mathcal{V}^*\Vec_\P1^{(r,k), \geq n}}|
  \mathcal{V}\in \Vec_\P1^{(r,d),\geq n}\}$ as
  irreducible components, each of dimension $-r^2$. Lemma \ref{IrrBun} follows.
\end{dem}
To get the theorem \ref{irr}, we have to describe concretely this
correspondence. Take $Z_1$ (resp. $Z_2$) an irreducible component of
$T^*\Vec_\P1^{(r,d-l)}$ (resp. $T^*\Coh_\P1^{(0,l)}$). It is the
closure of the conormal to strata $\mathcal{V}$ for $\mathcal{V}\in
\Vec_\P1^{(r,d-l)}$ by lemma \ref{IrrBun} (resp. $\CalO_\lambda$ for
$\lambda \vdash l$ by lemma \ref{lemmetorsion}). Hence $\pi_2^{-1}(Z_1\times Z_2)$ is an irreducible component of
$\underline{L}^{r,d,l}$ containing the substack whose objects are
$\{(\mathcal{V},\CalO_{x_1}^{(\lambda_1)}\oplus \cdots
\CalO_{x_{l(\lambda)}}^{(\lambda_{l(\lambda)})},f_1,f_2,f_3)|x_i \neq
    x_j\}$ as a dense substack. Then the substack $X_{\mathcal{V},\lambda}$
    of $\underline{\Lambda}_\P1^{r,d,l}$ is a dense substack of
    $\pi_1(\pi_2^{-1}(Z_1\times Z_2))$. This proves that the irreducible
    components of $\underline{\Lambda}_\P1^{r,d,l}$ are exactly
    $\overline{X}_{\mathcal{V},\lambda}\cap
    \underline{\Lambda}_\P1^{r,d,l}$ for $\mathcal{V}$ with
    $\mathcal{V}\in\Vec_\P1^{(r,d-l)}$ and $\lambda \vdash l$. We also
    have that $\underline{\Lambda}_\P1^{r,d,l}$ is pure of dimension
    $-r^2$. We have a locally finite union
\[
\underline{\Lambda}_\P1^{(r,d)}=\bigsqcup_{l\in \Nn}
\underline{\Lambda}_\P1^{r,d,l}
\]
which implies that $\underline{\Lambda}_\P1^{(r,d)}$ is pure of
dimension $-r^2$. As $\overline{X}_{\mathcal{V},\lambda}\cap
    \underline{\Lambda}_\P1^{r,d,l}$ is dense in
    $\overline{X}_{\mathcal{V},\lambda}$, the substack
    $\overline{X}_{\mathcal{V},\lambda}$ is an irreducible component
    of $\underline{\Lambda}_\P1^{(r,d)}$ and every irreducible
    component comes this way. Theorem \ref{irr} follows.

\section{The Higgs algebra}
\subsection{}
We consider the set $F(\Higgsa)$ of constructible functions on
$\Higgsa$, which are functions $f:\Higgsa \rightarrow \Qq$ satisfying:

- $\Im(f)$ is finite,

- $\forall c \in \Qq^*$, $f^{-1}(c)$ is constructible.\\
These objects are not true functions; such an object is
a partition of $\Higgsa$ in a finite number of constructible subsets
with a rational number attached to each element of the partition.\\
Define 
\[
F(\Higgs)=\bigoplus_{\alpha \in K^+(\Coh_\P1)} F(\Higgsa).
\]
We equip this space with a convolution product.\\
Consider the following diagram:
\[
\Higgsa \times \Higgsb  \xleftarrow{p_1} \text{\bf{E}}
\xrightarrow{p_2} \Higgsab
\]
where $\bf{E}$ parametrizes tuples 
\begin{align*}
\{(\CalF_1 \subseteq\CalF,f)|\CalF \in
\Coh_\P1,[\CalF]=\alpha+\beta,f \in
\Hom^\text{nilp}(\CalF,\CalF(-2)),\\ [\CalF_1]=\beta,f(\CalF_1)\subseteq
\CalF_1(-2)\}.
\end{align*}
The map $p_1$ is given by
$p_1(\CalF_1 \subseteq \CalF,f)=((\CalF_1,f_{|\CalF_1}),(\CalF/\CalF_1,f_{|\CalF/\CalF_1}))$.
  The map $p_2$ is $p_2(\CalF_1\subseteq \CalF,f)=(\CalF,f)$. It is a
  representable proper map.\\
Given $g_1$ in $F(\Higgsa)$ and $g_2$ in $F(\Higgsb)$, define the
product $g_1g_2$ by 
\[
g_1g_2(\CalF,f)=\int_{p_2^{-1}(\CalF,f)}g_2(\CalF_1,f_{\CalF_1})g_1(\CalF/\CalF_1,f_{\CalF/\CalF_1})
\]
where as usual for a constructible function $f$ on a space $E$ of
finite type, we
define $\int_E f=\sum_{a\in \Cc^*} a\chi(f^{-1}(a))$.\\
Now define a completion $\tilde{F}(\Higgs)$ of $F(\Higgs)$ as the inductive limit
\[
\tilde{F}(\Higgs)=\varinjlim_n F(\underline{\Lambda}_\P1^{\geq n})
\]
so that an element $g$ in $\tilde{F}(\Higgs)$ is a function whose restriction $g_{\geq n}$ to any open substack $\underline{\Lambda}_\P1^{\geq n}$ is a constructible function.\\
The convolution product in $F(\Higgs)$ then extends to
$\tilde{F}(\Higgs)$: for two functions $g_1,g_2$ we have
$g^1g^2(\CalF,f)=g^1_{\geq n}g^2_{\geq n}(\CalF,f)$ for $n<<0$
(depending on $\CalF,f$).\\
Note that after this completion the functions we consider are called
\textit{locally constructible} in \cite{Jo}.
\subsection{}
We have endowed the set $\tilde{F}(\Higgs)$ with a convolution
product. Define the Higgs algebra as the subalgebra $\alg$ generated
by the following elements:
\begin{enumerate}
\item $1_{(0,d)}=\chi_{\Coh_\P1^{(0,d)}}$ for $d \in \Nn^*$, the characteristic function of
    the zero section $\Coh_\P1^{(0,d)}\subseteq
    \underline{\Lambda}_\P1^{(0,d)}$ of the bundle
    $T^*\Coh_\P1^{(0,d)} \rightarrow \Coh_\P1^{(0,d)}$.
\item $1_{(1,n)}=\chi_{\Coh_\P1^{(1,n)}}$ for $n \in \Zz$, where
  $\Coh_\P1^{(1,n)}\subseteq \underline{\Lambda}_\P1^{(1,n)}$ is
  the zero section of the corresponding cotangent bundle.
\end{enumerate}
We have a natural decomposition:
\[
\alg=\bigoplus_{\alpha \in K^+(\Coh_\P1)}\alga
\]
corresponding to the decomposition $F(\Higgs)=\oplus_{\alpha \in K^+(\Coh_\P1)}F(\Higgsa)$.
\begin{prop}\label{rel}
These elements satisfy the following relations:

$(1)$ $1_{(0,d)}1_{(0,d')}=1_{(0,d')}1_{(0,d)}$ for every $d,d'\in
    \Nn^*$,

$(2)$ $1_{(0,d)}1_{(1,n)}=\sum_{k=0}^{n}(k+1)1_{(1,n+k)}1_{(0,d-k)}$
  for $d\in \Nn^*$ and $n\in \Zz$,

 \begin{multline*}(3) \ (n-l+2)(1_{(1,n-1)}1_{(1,l+1)}-1_{(1,l-1)}1_{(1,n+1)}) \\ =(n-l)(1_{(1,n)}1_{(1,l)}- 1_{(1,l-2)} 1_{(1,n+2)})
\end{multline*}
for $l,n\in \Zz$.
\end{prop}

\begin{dem}
Define $\underline{\Lambda}_{tor}$ to be the substack of $\Higgs$ parametrizing
pairs $(\CalF,f)$ with $\CalF$ a torsion sheaf.
We have the decomposition $F(\underline{\Lambda}_{tor})=\otimes_{x\in
  \P1}F(\underline{\Lambda}_{tor,x})$ where $\underline{\Lambda}_{tor,x}$ is the substack of
$\underline{\Lambda}_{tor}$ parametrizing couples $(\CalF,f)$ with
$\text{support}(\CalF)=x$.\\
The subalgebras $F(\underline{\Lambda}_{tor,x})$ commute with each other, so it remains to
study one such subalgebra.
If we want to compute the product $1_{(\CalF_1,f_1)}1_{(\CalF_2,f_2)}$
over an element $(\CalF,f)$, we have to compute the Euler
  characteristic of the set of subobjects $\CalF_2 \subseteq
  \CalF$ such that the following diagram commutes:
\[
\xymatrix{
0 \ar[r] & \CalF_2 \ar[r] \ar[d]^{f_2} & \CalF \ar[r] \ar[d]^f &
\CalF_1 \ar[r] \ar[d]^{f_1}& 0 \\
0 \ar[r] & \CalF_2 \ar[r] & \CalF \ar[r] &\CalF_1 \ar[r] & 0
}
\]
Apply the exact functor $\tilde{\_}=\Hom(\_,\Cc_x)$ where $\Cc_x$ is
the skyscraper sheaf at $x$ we have:
\[
\xymatrix{
0 \ar[r] & \tilde{\CalF}_1 \ar[r]  & \tilde{\CalF} \ar[r]  &
\tilde{\CalF}_2 \ar[r] & 0 \\
0 \ar[r] & \tilde{\CalF}_1 \ar[r] \ar[u]^{\tilde{f}_1} & \tilde{\CalF} \ar[r] \ar[u]^{\tilde{f}}&\tilde{\CalF}_2 \ar[r]\ar[u]^{\tilde{f}_2} & 0
}
\]
But as this functor preserves isomorphism classes, we have
$(\tilde{\CalF}_i,\tilde{f}_i)=(\CalF_i,f_i)$ and we see that the
number we are calculating is the same as in the product
$1_{(\CalF_2,f_2)}1_{(\CalF_1,f_1)}$.\\
Hence the product in $F(\underline{\Lambda}_{tor,x})$ is commutative, and hence so is the
product in $F(\underline{\Lambda}_{tor})$.\\
Now we prove the second relation.
We write the product:
\begin{align*}
1_{(0,d)}1_{(1,n)}=\sum_{(\CalF,f)}\chi(\CalF_2 \subseteq \CalF\
| \ \CalF_2 \text{ is of class } (1,n),& \\ &\hspace{-6cm}  \CalF/\CalF_2 \text{ is of class
  }(0,d), \  f_{|\CalF_2}=0, \ f_{|\CalF/\CalF_2}=0) 1_{(\CalF,f)}
\end{align*}
 The product $1_{(0,d)}1_{(1,n)}$ is
non-zero only on couples $(\CalF,f)$ with $\text{deg}(\CalF)=n+d$ and
$\text{rk}(\CalF)=2$. The coefficient of the product on the element of the
basis $1_{(\CalF,f)}$ is equal to the Euler characteristic of the set
of injections $\CalF_1 \hookrightarrow \CalF$, where $\CalF_1$ is a
sheaf of degree $n$ and of rank $1$, which make the following diagram
commutative:
\[
\xymatrix{
0 \ar[r] & \CalF_1 \ar[r] \ar[d]^{0} & \CalF \ar[r] \ar[d]^f &
\CalF_2 \ar[r] \ar[d]^{0}& 0 \\
0 \ar[r] & \CalF_1(-2) \ar[r] & \CalF(-2) \ar[r] &\CalF_2(-2) \ar[r] & 0
}
\]
We may rewrite the conditions imposed by the diagram as
\[
\Im(f)(2)\subseteq \CalF_1 \subseteq \Ker(f).
\]
 As $\CalF$ is of rank
one and $f$ is nilpotent, $\Im(f)$ is a torsion sheaf and
$\Im(f)(2)=\Im(f)$.\\
Define $\CalF^\text{fr}:=\CalF/\CalF^\text{tor}$.
The following lemma will often be used:
\begin{lem}\label{prodeuler}
For $\CalF$, $\CalG$ two coherent sheaves over a smooth projective
curve $X$, put:
\[
Gr_\CalF^\CalG:=\{\mathcal{H}\subseteq \CalG \ | \  \mathcal{H}\simeq
\CalF\} \ \ \text{  (a projective variety)}.
\]
Then $\chi(Gr_\CalF^\CalG)=\chi(Gr_{\CalF^{tor}}^{\CalG^{tor}})\chi(Gr_{\CalF^{fr}}^{\CalG^{fr}})$.
\end{lem}
\begin{dem}
Fix a decomposition:$\begin{cases}
 \CalF=\CalF^\text{fr} \oplus
\CalF^\text{tor}\\
\CalG=\CalG^\text{fr} \oplus \CalG^\text{tor}
\end{cases}$.\\
We have:
\[
\Hom(\CalF,\CalG)=\Hom(\CalF^\text{fr},\CalG^\text{fr})\oplus
\Hom(\CalF^\text{tor},\CalG^\text{tor}) \oplus
\Hom(\CalF^\text{fr},\CalG^\text{tor})
\]
\[
\Aut(\CalF)=(\Aut(\CalF^\text{fr}) \times
\Aut(\CalF^\text{tor}))\rtimes \Hom(\CalF^\text{fr},\CalF^\text{tor}).
\]
Write $\Hom(\CalF,\CalG)^\text{inj}$ for the subset of injections.
We have
$\Hom(\CalF,\CalG)^\text{inj}=\Hom(\CalF^\text{fr},\CalG^\text{fr})^\text{inj}\times
\Hom(\CalF^\text{tor},\CalG^\text{tor})^\text{inj} \times
\Hom(\CalF^\text{fr},\CalG^\text{tor})$.\\
Since $\Aut(\CalF)$ acts freely on $\Hom(\CalF,\CalG)^\text{inj}$, we
have:
\begin{align*}
\chi(Gr_\CalF^\CalG) & =\chi(\Hom(\CalF^\text{fr},\CalG^\text{fr})^\text{inj}\times
\Hom(\CalF^\text{tor},\CalG^\text{tor})^\text{inj} \times
\Hom(\CalF^\text{fr},\CalG^\text{tor})/\chi(\Aut(\CalF)) \\
 & =\chi(\Hom(\CalF^\text{fr},\CalG^\text{fr})^\text{inj})/\chi(\Aut(\CalF^\text{fr}))
.
\chi(\Hom(\CalF^\text{tor},\CalG^\text{tor})^\text{inj})/\chi(\Aut(\CalF^\text{tor})\\
 & =\chi(Gr_{\CalF^{tor}}^{\CalG^{tor}})\chi(Gr_{\CalF^{fr}}^{\CalG^{fr}}).
\end{align*}
\end{dem}
Using lemma \ref{prodeuler}, we have:
\begin{multline*}
\chi(\{\Im(f)(2)\subseteq \CalF_1 \subseteq
\Ker(f)\}) \\=\chi(\{\Im(f)(2)\subseteq \CalF_1^\text{tor} \subseteq
\Ker(f)^\text{tor}\}). \chi(\{\CalF_1^\text{fr} \subseteq \Ker(f)^\text{fr}\}).
\end{multline*}
We define $n_1=\text{deg}(\CalF_1^\text{tor})$ and
$n_2=\text{deg}(\CalF_1^\text{fr})$.\\
In the same way,  $k_1=\text{deg}(\Ker(f)^\text{tor})$,
$k_2=\text{deg}(\Ker(f)^\text{fr})$ and $t=\text{deg}(\Im(f))$.\\
We have the following: 
\[
n_1+n_2=n, \  k_1 +k_2+t=d+n, \  n_2 \leq k_2, \  t\leq n_1\leq k_1.
\]
If we denote by $Gr_{k}^\CalG$ the
Grassmanian of subsheaves of a torsion sheaf $\CalG$ of a given
degree $k$, we have that
\[
\begin{array}{ll}
\chi(\{\Im(f)(2)\subseteq \CalF_1^\text{tor} \subseteq
\Ker(f)^\text{tor}\}) & =\chi(\{ \CalF_1^\text{tor} \subseteq
\Ker(f)^\text{tor}/\Im(f)\})
\\ & =\chi(Gr_{n_1-t}^{\Ker(f)^{\text{tor}}/\Im(f)}).
\end{array}
\]
The subsheaf $\CalF_1^\text{fr}$ is just $\CalO(n_2)$. So
$Gr_{\CalF_1^\text{fr}}^{\CalF^\text{fr}}$ is
\[
(\Hom(\CalO(n_2),\CalF^{fr})-0)/\Cc^*=\mathbb{P}^{k_2-n_2}
\]
whose Euler characteristic is $k_2-n_2+1$ (if $k_2 \geq n_2$).\\
The product can finally be written:
\begin{multline}\label{torfr}
1_{(0,d)}1_{(1,n)} =\sum_{\CalF,f}\sum_{n_1=t}^{k_1}\chi(\mathbb{P}^{k_2-n_2})\chi(Gr_{n_1-t}^{\Ker(f)^{\text{tor}}/\Im(f)})1_{(\CalF,f)}
\\ = \sum_{\CalF,f}\sum_{n_2=k_2-d}^{k_2}(k_2-n_2+1)\chi(Gr_{n_1-t}^{\Ker(f)^{\text{tor}}/\Im(f)})1_{(\CalF,f)}
\end{multline}
where the sum is over the couples $(\CalF,f)$ with $\CalF$ of degree
$n+d$ and rank $1$, and $f\in \Hom(\CalF,\CalF(-2))$ nilpotent such
that $f^2=0$.\\
Let us now compute the product $1_{(1,n)}1_{(0,d)}$. The computation is
the same, except for the fact that now the subsheaf $\CalF_1$ is torsion
of degree $d$. We have to count the Euler characteristic of the set of
subobjects $\CalF_1$ of $\CalF$ such that $\Im(f) \subseteq \CalF_1
\subseteq \Ker(f)$. We get
\[
1_{(1,n)}1_{(0,d)}=\sum_{\CalF,f}\chi(Gr_{d-t}^{\Ker(f)^{\text{tor}}/\Im(f)})1_{(\CalF,f)}
\]
where the sum is over couples $(\CalF,f)$ as before with the
additional conditions $\text{degree}(\CalF^\text{tor})\geq d$ and
$\text{degree}(\Im(f))\leq d$.
To prove the second formula of the proposition, we rewrite
equation (\ref{torfr}). Introduce $C=k_2-n_2$ and exchange the sums as follows:
\begin{multline*}
1_{(0,d)}1_{(1,n)}=\sum_{\CalF,f}\sum_{C=d-k_1}^{d-t}(C+1)\chi(Gr_{C-k_2+n}^{\Ker(f)^{\text{tor}}/\Im(f)})1_{(\CalF,f)}\\
  =\sum_{C=0}^{d}(C+1)\sum_{\substack{\CalF\\d(\CalF^\text{tor})\geq
  d-c}}\sum_{f,\  d(\Im(f))=t\leq d-c}
  \chi(Gr_{C-k_2+n}^{\Ker(f)^\text{tor}/\Im(f)})1_{(\CalF,f)} 
\end{multline*}
We recognize the product $1_{(1,n+C)}1_{(0,d-C)}$ in the sum, so
we finally have:
\[
1_{(0,d)}1_{(1,n)}=\sum_{C=0}^{d}(C+1)1_{(1,n+C)}1_{(0,d-C)}
\]
as wanted.\\
Now we prove the third relation.
Let us determine the coefficient of the product $1_{(1,n)}1_{(1,l)}$ on
the element $1_{(\CalF,f)}$. For this coefficient to be non zero, the
sheaf $\CalF$ has to be of class $(2,l+n)$ and $f$ such that $f^2=0$. The
rank of $\Im(f)$ is at most one. We split the argument in two cases:
\\

\underline{Case 1}: rk($\Im(f)$)=1.\\
Using lemma (\ref{prodeuler}), we split the condition $\Im(f)(2)
\subseteq \CalF_1 \subseteq \CalF$ in two conditions:
\[
\Im(f)^\text{tor} \subseteq \CalF_1^\text{tor} \subseteq \Ker(f)^\text{tor}
\]
and
\[
\Im(f)^\text{fr}(2) \subseteq \CalF_1^\text{fr} \subseteq \Ker(f)^\text{fr}
\]
The coefficient we seek is the product of the Euler characteristic of
the set of 
$\{\CalF_1^\text{tor}\}$ and $\{\CalF_1^\text{fr}\}$ verifying the conditions
above.
The Euler characteristics are:
\[
\chi(Gr_{k-t}^{\Ker(f)^\text{tor}/\Im(f)^\text{tor}})
\]
and
\[
\chi(Gr_{l-k-(t'+2)}^{\Ker(f)^\text{fr}/\Im(f)^\text{fr}(2)})
\]
where $t'=d(\Im(f)^\text{fr})$.\\
Then we have
\[
1_{(1,n)}1_{(1,l)}(\CalF,f)= \sum_{k=t}^{d(\Ker (f)^\text{tor})}\chi(Gr_{k-t}^{\Ker(f)^\text{tor}/\Im(f)^\text{tor}}) \chi(Gr_{l-k-(t'+2)}^{\Ker(f)^\text{fr}/\Im(f)^\text{fr}(2)}).
\]
Now we use the duality for torsion sheaves: for any torsion
sheaf $\CalG$ we have $\chi(Gr_m^{\CalG})=\chi(Gr_{d(\CalG)-m}^{\CalG})$.\\
We apply this formula to the first coefficient in the right-hand side:
\[
\chi(Gr_{k-t}^{\Ker(f)^\text{tor}/\Im(f)^\text{tor}})=\chi(Gr_{k_1-k}^{\Ker(f)^\text{tor}/\Im(f)^\text{tor}})
\]
(we define as before $k_1=d(\Ker(f)^\text{tor})$)\\
For the second coefficient we have:
\[
\chi(Gr_{l-k-(t'+2)}^{\Ker(f)^\text{fr}/\Im(f)^\text{fr}(2)})=\chi(Gr_{k_2+k-l}^{\Ker(f)^\text{fr}/\Im(f)^\text{fr}(2)})
\]
Now we change $k$ into $k'=k_1-k+t$, and we use the equality
$k_1+k_2+t+t'=l+n$ to obtain the followoing expression in the
right-hand side:
\[
\chi(Gr_{k-t}^{\Ker(f)^\text{tor}/\Im(f)^\text{tor}})
\chi(Gr_{n-k-t'}^{\Ker(f)^\text{fr}/\Im(f)^\text{fr}(2)})
\]
We see that this is exactly the coefficient obtained when we compute
the product $1_{(1,l-2)}1_{(1,n+2)}(\CalF,f)$.\\
So if $\text{rk}(\Im(f))=1$ then
\[
(1_{(1,n)}1_{(1,l)}-1_{(1,l-2)}1_{(1,n+2)})(\CalF,f)=0.
\]

\underline{Case 2}: rk($\Im(f)$)=0.\\
In this case the conditions are:
\[
\Im(f) \subseteq \CalF_1^\text{tor} \subseteq \Ker(f)^\text{tor}
\]
and
\[
\CalF_1^\text{fr} \subseteq \Ker(f)^\text{fr}.
\]
Set $k=d(\CalF_1^\text{tor})$, which goes from $0$ to
$d(\Ker(f)^\text{tor})$, and $t=d(\Im(f)^\text{tor})=d(\Im(f))$. We
have for the Euler characteristic of the first set:
\[
\chi(\{\CalF_1^\text{tor}\})=\chi(Gr_{k-t}^{\Ker(f)^\text{tor}/\Im(f)}).
\]
For the second one, if we write $\Ker(f)=\CalO(v_1)\oplus
\CalO(v_2)\oplus \tau$, where $\tau$ is a torsion sheaf of degree $d$
and $v_1 \leq v_2$ then
\[
\chi(\{ \CalF_1^\text{fr} \})= \begin{cases}
  v_1+v_2-2(l-k)+2 & \text{if $l-k \leq v_1+1$} \\
  v_2-(l-k)+1 & \text{if $v_1<l-k$ and $l-k \leq v_2+1$}\\
  0 & \text{if $l-k>v_2$}
\end{cases}
\]
Therefore
\[
1_{(1,n)}1_{(1,l)}(\CalF,f)=\sum_{k=t}^{d(\Ker (f)^\text{tor})}S_k(\CalF,f)\chi(Gr_{k-t}^{\Ker(f)^\text{tor}/\Im(f)})
\]
where $S_k(\CalF,f):=\chi(\{ \CalF_1^\text{fr} \})$ is as above.\\
Now we compute $1_{(1,l-2)}1_{(1,n+2)}(\CalF,f)$. Arguing the same way,
\[
1_{(1,l-2)}1_{(1,n+2)}(\CalF,f)=\sum_{k=t}^{d(\Ker (f)^\text{tor})}S_k'(\CalF,f)\chi(Gr_{k-t}^{\Ker(f)^\text{tor}/\Im(f)})
\]
where $S_k'(\CalF,f)$ is given by:
\[
S_k'= \begin{cases}
  v_1+v_2-2(n+2-k)+2 & \text{if $n+2-k \leq v_1+1$} \\
  v_2-(l-k)+1 & \text{if $v_1<n+2-k$ and $n+2-k \leq v_2+1$}\\
  0 & \text{if $n+2-k>v_2$}
\end{cases}
\]
We set $m=d(\Ker f^\text{tor})+t-k$.\\
Note that
$\chi(Gr_{k-t}^{\Ker(f)^\text{tor}/\Im(f)})=\chi(Gr_{m-t}^{\Ker(f)^\text{tor}/\Im(f)})$. Thus
\[
1_{(1,l-2)}1_{(1,n+2)}(\CalF,f)=\sum_{m=t}^{d(\Ker (f)^\text{tor})}S_{d(\Ker f^\text{tor})+t-m}'(\CalF,f)\chi(Gr_{m-t}^{\Ker(f)^\text{tor}/\Im(f)}).
\]
Using $v_1+v_2+k+m=v_1+v_2+d(\Ker f^\text{tor})+t=n+l$. We may rewrite the conditions for the values of $S_{d(\Ker f^\text{tor})+t-m}'$ in terms of $m$:
\[
S_{d(\Ker f^\text{tor})+t-m}'= \begin{cases}
  v_1+v_2-2(n-k)-2 & \text{if $m-l > v_2$} \\
  v_2-(n-k)-1 & \text{if $l-m\leq v_2$ and $l-m>v_1$}\\
  0 & \text{if $l-m \leq v_1$}
\end{cases}
\]
Hence $S_k-S_{d(\Ker f^\text{tor})+t-k}'=2(n-l+2)$ for any $k$ and 
\[
(1_{(1,n)}1_{(1,l)}-1_{(1,l-2)}1_{(1,n+2)})(\CalF,f)=2(n-l+2)\sum_{k=t}^{d(\Ker (f)^\text{tor})}\chi(Gr_{k-t}^{\Ker(f)^\text{tor}/\Im(f)})
\]
We deduce
\begin{multline*}
1_{(1,n)}1_{(1,l)}-1_{(1,l-2)}1_{(1,n+2)}= \\ 2\sum_{\substack{\CalF,f \\ rk(\Im(f))=0}}
\left( \sum_{k=t}^{d(\Ker(f)^\text{tor})}
\chi(Gr_{k-t}^{\Ker(f)^\text{tor}/\Im(f)})(n-l+2)\right) 1_{(\CalF,f)}
\end{multline*}
If we set
\[
\phi_p=2\sum_{\substack{\CalF,f \\ rk(\Im(f))=0}}
 \sum_{k=t}^{d(\Ker(f)^\text{tor})}
\chi(Gr_{k-t}^{\Ker(f)^\text{tor}/\Im(f)}) 1_{(\CalF,f)}
\]
where the sum is over pairs $(\CalF,f)$ where $\CalF$ of rank 2 and
degree p, and $f^2=0$, we can write:
\[
1_{(1,n)}1_{(1,l)}-1_{(1,l-2)}1_{(1,n+2)}=(n-l+2)\phi_{n+l}
\]
In the same manner:
\[
1_{(1,n-1)}1_{(1,l+1)}-1_{(1,l-1)}1_{(1,n+1)}=(n-l)\phi_{n+l}
\]
The third relation follows.

\end{dem}

\subsection{}
Using the relation (3) of Proposition \ref{rel}, we see that if $n+l=2k$, with
$k\in \Zz$, then
\[
1_{(1,n)}1_{(1,l)}=1_{(1,l-2)}1_{(1,n+2)}+\frac{n-l+2}{2}(1_{(1,k)}^2-1_{(1,k-2)}1_{(1,k+2)})
\]
and if $n+l=2k+1$, with $k\in\Zz$, we have:
\[
1_{(1,n-1)}1_{(1,l+1)}=1_{(1,l-1)}1_{(1,n+1)}+(n-l+2)(1_{(1,k)}1_{(1,k+1)}-1_{(1,k-1)}1_{(1,k+2)}).
\]
So it is always possible to express a product of $1_{(1,n_i)}$ as a
linear combination of ordered products $\prod 1_{(1,l_i)}$ with
$l_1\leq l_2\leq\cdots$.\\
Combined with the relation (1) of Prop \ref{rel}, we deduce that
the multiplication induces a surjective morphism of vector spaces:
\begin{equation*}\label{ord}
\mathcal{H}^{fr}\otimes \mathcal{H}^{tor} \twoheadrightarrow \alg
\end{equation*}
where $\mathcal{H}^{fr}$ is the space 
\[
\mathcal{H}^{fr}= \sum_{l_1 \geq l_2 \geq \cdots \geq
  l_s}\Cc1_{(1,l_1)}\cdots 1_{(1,l_s)}
\] and
$\mathcal{H}^{tor}$ is the subalgebra generated by the elements
$1_{(0,d)}$.

\section{The semicanonical basis}\label{semicanbasis}
We want to construct a basis of our algebra indexed by the irreducible
components of $\Higgs$ in the spirit of \cite{L1}.\\
More precisely if $h \in F(\Higgsa)$, for every $Z \in
\text{Irr }\Higgsa$ there is a unique $c\in \Qq$ such that
$h^{-1}(c)\cap Z$ is an open dense subset of $Z$. We define
$\rho_Z(h)=c$.\\
In the following, we denote by $Z_{(\underline{n},\lambda)}$ the
irreducible component $\overline{X}_{\mathcal{V},\lambda}$, where
$\lambda$ is a partition and $\underline{n}$ is a $\Zz$-partition
corresponding to the vector bundle $\mathcal{V}=\oplus_i \CalO(n_i)$.
\begin{thm}\label{irred}
For each irreducible component $Z=Z_{(\underline{n},\lambda)}\in
\Lambda^\alpha$ there exists an element $f_{(\underline{n},\lambda)}=f_Z \in
\alga$
with $\rho_Z(f_Z)=1$ and $\rho_{Z'}(f_Z)=0$ for every $Z' \neq Z$.
\end{thm}
We fix some $\alpha=(r,d)$.\\
Assume first that $r=0$, hence the irreducible components are parametrized by
partitions of $d$. Define for a partition $\lambda$ of $d$:
\[
F_\lambda=\{ (\CalO_{x_1}^{(\lambda_1)}\oplus \cdots \oplus
\CalO_{x_d}^{(\lambda_d)}, f)\in \Lambda^{(0,d)}| x_i \neq x_j, f
\text{ generic}\}
\]
(here ``generic'' means that the degree of $\Im f$ is maximal) 
which is an open dense subset of the irreducible component associated
to $\lambda$. We want to find for each $\lambda \vdash d$ an element
$h \in \alg^{(0,d)}$ such that $h$ is 1 on $F_\lambda$ and 0 on any other
$F_\mu$, $\mu \neq \lambda$.\\
Define $\lambda'$ to be the transpose of $\lambda$ and $l(\lambda)$
the length of $\lambda$. Define also
$1_\lambda=\prod_{i=1}^{l(d)} 1_{(0,\lambda_i)}\in \alg^{(0,d)}$, denote
by $\lambda \prec \nu$ the Bruhat order on partitions of $d$ and by
$K$ the set of functions on $\Lambda^{(0,d)}$ which are generically 0.
\begin{lem}
For every $\lambda \vdash d$ we have:
\[
1_{\lambda'}\in 1_{F_\lambda} +\oplus_{\mu
  \prec\lambda}\Cc1_{F_\mu} + K 
\]
\end{lem}
\begin{dem}
On a point $(\CalF,f)\in F_\mu$, the type of the nilpotent operator
$f$ acting on $H^0(\CalF)$ is $\mu$. Note that $\dim H^0(\CalF)=\text{deg}(\CalF)=d$. We are interested in the set of $\mu$s such that $1_{\lambda'}$
is non zero on $F_\mu$.
For such a $\mu$ and $(\CalF,f) \in F_\mu$ there exists a filtration 
\[
V_0\subseteq V_1
\subseteq \cdots \subseteq V_{l(\lambda')}=V:=H^0(\CalF)
\]
 with $f(V_i)\subseteq V_{i-1}$ and
$\dim V_i/V_{i-1}=\lambda_i'$ for $1 \leq i\leq l(\lambda')$. So we have for
every $i$, $\dim\Ker f^i\geq \lambda_1'+\cdots +\lambda_i'$. As $f$ is
of type $\mu$, this means that $\lambda'\preccurlyeq\mu'$, hence $\mu
\preccurlyeq \lambda$. In the case of equality, the filtration is
uniquely determined by $f$. Thus $1_{\lambda'}$ is 1 on $F_\lambda$
and is non zero on $F_\mu$ (with $\mu \neq \lambda$) only if $\mu
\prec \lambda$. 
\end{dem}
To obtain an element $h$ which is generically 1 on $F_\lambda$ and 0
on $F_\mu$, $\mu \neq \lambda$, we proceed by induction on $\lambda$
with the Bruhat order by using the preceding lemma.\\
This settles the case $r=0$.\\
The proof for the case $r>0$ splits into several steps. We begin with a lemma, in which we are
only interested in irreducible components of the form
$Z_{(\underline{n},\lambda)}$ with $\lambda=0$:
\begin{lem}\label{vect}
For every $Z=Z_{(\underline{n},0)}$ and $k \in \Zz$ there exists
an element $g_{Z,k} \in \alga$ such that:
\begin{itemize}
\item $\rho_Z(g_{Z,k})=1$ if $Z \cap
  \Lambda_\alpha^{\geq k}$ is dense in $Z$.
\item $\rho_{Z'}(g_{Z,k})=0$ if
  $Z'=Z_{(\underline{n}',0)}\neq Z$ and $Z' \cap
  \Lambda_\alpha^{\geq k}$ is dense in $Z'$.
\item $g_{Z,k-1}|_{\Lambda_\alpha^{\geq k}}=g_{Z,k}|_{\Lambda_\alpha^{\geq k}}$.
\end{itemize}
\end{lem}
\begin{dem}
We introduce $1_{\underline{n}}=1_{(1,n_1)}\cdots 1_{(1,n_r)}$. Define
an order on $\Zz$-valued (anti-)partitions of $d$ as follows:
$\underline{n} \prec \underline{m}$ if there exists $j\geq 1$ such
that $n_i=m_i$ for $i>j$ and $n_j<m_j$. The construction of our
function $g_{Z,k}$ will follow from the next lemma:
\begin{lem}
If the function $1_{\underline{n}}$ is generically non zero on
$Z_{(\underline{m},0)}$ then $\underline{m}=\underline{n}$ or
$\underline{n}\prec \underline{m}$. Moreover it is generically non zero on $Z_{(\underline{n},0)}$.
\end{lem}
\begin{dem}
To see this, remark that if $1_{\underline{n}}$ is non zero on
$(\CalF,f)\in Z_{(\underline{m},0)}$ (with $\CalF$ a vector bundle)
then as it appears in the product
$(1_{(1,n_1)} \cdots 1_{(1,n_{r-1})})1_{(1,n_r)}$, we have an
injection $\CalO(n_r)\hookrightarrow \CalF$. So we have $n_r \leq
m_r$.\\
Two cases appear:\\
First case: $n_r<m_r$. Then $\underline{n}\prec \underline{m}$.\\
Second case: $n_r=m_r$. In this case the quotient $\CalF/\CalO(n_r)$
is still a vector bundle, coming from the product
$1_{(1,n_1)}\cdots 1_{(1,n_{r-1})}$. An easy induction on the rank $r$
gives the first result.\\
To see the second part, as the set $\{\CalO(n_1)\oplus \cdots \oplus
\CalO(n_r),f\}$ is dense in $Z_{(\underline{n},0)}$, we only have to
compute $1_{\underline{n}}(\CalF,f)$ for $\CalF=\CalO(n_1) \oplus
\cdots \oplus \CalO(n_r)$. We are in the second case:
$1_{\underline{n}}(\CalF,f)$ is equal to the Euler caracteristic of the
set of subobjects $\CalO(n_r) \subseteq \CalF$ (which is non
zero integer) multiplied by
$1_{(1,n_1)}\cdots 1_{(1,n_{r-1})}(\CalF/\CalO(n_r),f')$. An easy
induction on $r$ gives the result.
\end{dem}
Now the lemma \ref{vect} comes from an induction on the order $\prec$,
or more precisely on $d(\underline{m},\underline{n})$, where we define
for $\underline{n} \prec \underline{m}$ the integer
$d(\underline{m},\underline{n})$ as the maximal length $c$ of a chain
$\underline{n}=\underline{n}_0 \prec \cdots \prec
\underline{n}_c=\underline{m}$. We define a sequence:
\begin{itemize}
\item $g_{Z,k}^{(0)}=\rho_Z(1_{\underline{n}})^{-1}1_{\underline{n}}$
\item \begin{displaymath}g_{Z,k}^{(j)}=g_{Z,k}^{(j-1)}
  -\sum_{d(\underline{m},\underline{n})=j}\rho_{Z_{(\underline{m},0)}}(g_{Z,k}^{(j-1)})1_{\underline{m}}.
\end{displaymath}
\end{itemize}
By construction, this sequence has the property:
$g_{Z,k}^{(j)}$ is generically $1$ on $Z$ and is
generically non zero on $Z_{(\underline{m},0)}$ only if
$\underline{n}\prec \underline{m}$ and
$d(\underline{m},\underline{n})>j$.\\
As we consider irreducible components $Z'=Z_{(\underline{m},0)}$ such
that $Z' \cap \Lambda_\alpha^{\geq k}$ is dense in $Z'$, the set of
integers $d(\underline{m},\underline{n})$ is finite (as the $m_i$'s
are bounded below and the sum is equal to $d$) and we denote by $h$
its maximal value. Define $g_{Z,k}=g_{Z,k}^{(h)}$. It is clearly a
solution to our problem. The last property in the lemma is also
clearly verified.
\end{dem}
We now give a refinement of the preceding lemma:
\begin{lem}\label{vect2}
For every $Z=Z_{(\underline{n},\lambda)}$ and $k \in \Zz$ there exists
an element $g_{Z,k}'$ such that:
\begin{itemize}
\item $\rho_Z(g_{Z,k}')=1$ if $Z \cap
  \Lambda_\alpha^{\geq k}$ is dense in $Z$.
\item $\rho_{Z'}(g_{Z,k}')=0$ if
  $Z'=Z_{(\underline{n}',\lambda')}$ is such that $|\lambda'|\leq
  |\lambda|$, $Z'\neq Z$ and $Z' \cap
  \Lambda_\alpha^{\geq k}$ is dense in $Z'$.
\item $g_{Z,k-1}'|_{\Lambda_\alpha^{\geq
      k}}=g_{Z,k}'|_{\Lambda_\alpha^{\geq k}}$.
\end{itemize}
\end{lem}
\begin{dem}
Set $g_{Z,k}'=g_{Z_{(\underline{n},0)}}1_\lambda$, where
$1_\lambda$ has been defined in the proof of Theorem \ref{irred} for torsion sheaves. Take a point $(\CalF,f)$ in the generic part of
$Z'=Z_{(\underline{n}',\lambda')}$ and assume that $g_{Z,k}'(\CalF,f)\neq
0$. By definition there is an injection $\tau \hookrightarrow \CalF^{tor}$,
where $\tau$ is a torsion sheaf such that $(\tau,f|_\tau)$ is in the
support of $1_\lambda$. As we take $(\CalF,f)$ in the generic part of
$Z'$, the element $(\tau, f|_\tau)$ is in the generic part of
$Z_{(0,\lambda)}$, and $1_\lambda$ is generically non zero only on
$Z_{(0,\lambda)}$. So we may assume that $\tau=\tau_\lambda$.\\
We have two cases:\\
The injection $\tau_\lambda \hookrightarrow \CalF^{tor}$ is not an
isomorphism. Then $|\CalF^{tor}|>|\lambda|$ and $\lambda'>\lambda$.\\
The injection $\tau_\lambda \hookrightarrow \CalF^{tor}$ is an
isomorphism. We have $\lambda'=\lambda$. The quotient
$\CalF'=\CalF/\tau_\lambda$ is a vector bundle such that $(\CalF',f')$
is in the generic support of $g_{Z_{(\underline{n},,0)}}$, which implies
by lemma \ref{vect} that $\underline{n}'=\underline{n}$ and the result
follows. We have of course that $\rho_Z(g_{Z,k}')=1$.
\end{dem}
We can now prove Theorem \ref{irred}.\\
First remark that on $\Lambda_\alpha^{\geq k}$ the degree of the
torsion part of a sheaf is bounded by $d-rk$. We denote this number by
$\tau(k)$. Recall the fixed component $Z=Z_{(\underline{n},\lambda)}$. We construct a
sequence of functions $f_{Z,k}^{(j)}$ as follows:
\begin{itemize}
\item $f_{Z,k}^{(0)}=g_{Z,k}'$
\item \begin{displaymath}f_{Z,k}^{(j)}=f_{Z,k}^{(j-1)}-
  \sum_{|\mu|=|\lambda|+j}\rho_{Z_{(m,\mu)}}(f_{Z,k}^{(j-1)})g_{(Z_{(m,\mu)},k)}' \end{displaymath}
\end{itemize}
We define $f_{Z,k}=f_{Z,k}^{(\tau(k)-|\lambda|)}$.\\
By construction we see that $f_{Z,k}^{(j)}$ has the following
properties:
\begin{enumerate}
\item $f_{Z,k}^{(j)}$ is generically $1$ on
  $Z=Z_{(\underline{n},\lambda)}$
\item $f_{Z,k}^{(j)}$ is generically $0$ on each
  $Z'=Z_{(\underline{m},\mu)}\neq Z$ such that $|\mu|\leq |\lambda|+j$.
\end{enumerate}
Then by definition the function $f_{Z,k}=f_{Z,k}^{(\tau(k)-|\lambda|)}$ is
generically $1$ on $Z$ and $0$ on each $Z_{(\underline{m},\mu)}$ such
that $|\mu|\leq \tau(k)$, and so on each $Z'$ such that $Z'\cap
\Lambda_\alpha^{\geq k}$ is dense in $Z'$.\\
It remains to verify that $f_{Z,k}$ coincides with $f_{Z,k-1}$ on
$\Lambda_\alpha^{\geq k}$ to go to the limit. But the lemma \ref{vect2}
gives us that $f_{Z,k}^{(0)}$ and $f_{Z,k-1}^{(0)}$ coincides on
$\Lambda_\alpha^{\geq k}$. By construction of our sequences, we see
that for $j\leq \tau(k)-|\lambda|$, the fact that $f_{Z,k}^{(j)}$ and
$f_{Z,k-1}^{(j)}$ coincide on $\Lambda_\alpha^{\geq k}$ implies that
$f_{Z,k}^{(j+1)}$ and $f_{Z,k-1}^{(j+1)}$ coincide on
$\Lambda_\alpha^{\geq k}$. Then for $j>\tau(k)-|\lambda|$ the
  functions $f_{Z,k-1}^{(j)}$ and $f_{Z,k-1}^{(j+1)}$ clearly coincide
  on $\Lambda_\alpha^{\geq k}$.\\
So we have that $f_{Z,k-1}|_{\Lambda_\alpha^{\geq
    k}}=f_{Z,k}|_{\Lambda_\alpha^{\geq k}}$, then the limit
\[
f_Z=\lim_{k \rightarrow -\infty}f_{Z,k}
\]
is well defined, belongs to $\alga$ and clearly verifies the
properties in the theorem.

\section{The affine Lie algebra $\widehat{sl}_2$}

In this section we recall the definition of the affine Lie algebra
$\hat{sl}_2$, its positive part and the completion of the positive
part.

We define the Lie algebra $\widehat{sl}_2$ as the loop algebra of $sl_2$
(rather than its Kac-Moody definition). This is the Lie algebra $sl_2
\otimes \Cc[t,t^{-1}]$ with Lie bracket given by $[a \otimes t^i,b
\otimes t^j]=[a,b]\otimes t^{i+j}$.

Its roots lattice is $\widehat{Q}=\Zz\alpha_1 \oplus \Zz\delta$, where
$\alpha_1$ is the simple real positive root and $\delta$ is the
indivisible imaginary root. The Drinfeld (non-standard) set of
positive roots is $\widehat{Q}^+=\{\alpha_1+\Zz\delta\}\cup \{
\Nn^*\delta\}$.\\
We have an isomorphism $\rho : K_0(\P1) \rightarrow \widehat{Q}$ given by 
\[
\rho([\CalF])=\text{rk}(\CalF)\alpha_1+\text{deg}(\CalF)\delta.
\]
This isomorphism restricts to an isomorphism between positive parts:
\[
\rho: K_0^+(\P1) \simeq \widehat{Q}^+
\]
The Drinfeld positive part $\textbf{U}^+$ of the enveloping Lie algebra $\textbf{U}(\widehat{sl}_2)$
is the subalgebra generated by $e_i:=e\otimes t^i$, $i\in\Zz$ and
$h_j:=h\otimes t^j$, $j\in \Nn^*$.\\
The enveloping algebra $\textbf{U}(\widehat{sl}_2)$ is $K_0(\P1)\simeq
\widehat{Q}$-graded, and the (Drinfeld) positive part is
$K_0^+(\P1)$-graded.\\
For $\alpha= k\alpha_1+l\delta \in \widehat{Q}^+$, define $deg(\alpha)=l$.\\
Consider a completion $\widehat{\bf U}^+(\widehat{sl}_2)$ of
$\textbf{U}^+(\widehat{sl}_2)$ as follows. We write
$\widehat{\bf U}^+(\widehat{sl}_2)=\oplus_{\gamma \in K_0^+(\P1)}
\widehat{\bf U}^+(\widehat{sl}_2)[\gamma]$ where
\begin{align*}
\widehat{\bf U}^+(\widehat{sl}_2)[\gamma]=\{ \sum_{i \geq 0}a_ib_i \ | \ a_i
\in \textbf{U}^+(\widehat{sl}_2)[\alpha_i], \ b_i
\in \textbf{U}^+(\widehat{sl}_2)[\beta_i], \\ \alpha_i+\beta_i=\gamma, \  deg(\beta_i) \rightarrow \infty \}.
\end{align*}
The product is well defined in this completion.\\
For a positive integer $l$ define $P_l(X_1,\cdots, X_l)$ by the
formula: $P_l(X_1,\cdots, X_l)=\sum_{\lambda \vdash
  l}\prod_{i}X_{\lambda_i}$.\\
We can now state our first main theorem:
\begin{thm}
(1) The map
\[
\begin{array}{llll}
\psi: & \alg_\P1 &\rightarrow &\widehat{\bf{U}}^+(\widehat{sl}_2)\\
& 1_{(0,d)}&\mapsto &P_d(h_1, \cdots,h_d) \\
& 1_{(1,n)} &\mapsto &\sum_{m\geq 0} e_{n-m}P_m(h_1, \cdots, h_m)
\end{array}
\]
extends to an isomorphism of algebras.

(2) There exists a unique basis $B=\{b_Z\}$ of $\alg_\P1$, called the
    semicanonical basis, which is parametrised by irreducible
    components $Z$ of $\Lambda$, such that:
\begin{itemize}
\item $b_Z$ is generically $1$ on $Z$.
\item $b_Z$ is generically $0$ on $Z' \neq Z$.
\end{itemize}
\end{thm}
Let $\textbf{H}_\P1$ be the Hall algebra of locally constructible
functions on $\Coh_\P1$ (see section \ref{semicanbasis}). We use the well known isomorphism for $\textbf{H}_\P1$.\\
The next theorem is a version of \cite{Sc2}, theorem 10.10 for $q=1$,
itself a variant of \cite{Kap} and \cite{BK}.
\begin{thm}\cite{BK}
The assignement $e_l \mapsto [\CalO(l)]$ for $l \in \Zz$, $h_r \mapsto
T_r$ extends to an algebra isomorphism:
\[
\phi: \widehat{\bf U}^+(\widehat{sl}_2) \rightarrow \bf H_{\P1}
\]
\end{thm}
Now we consider the map
\[
H: \alg_\P1 \rightarrow \bf H_{\P1}
\]
defined by $H(g)(\CalF)=g(\CalF,0)$.

\begin{lem}
We have the following:
\begin{enumerate}
\item The map $H$ is a morphism of algebras.
\item $H$ is surjective.
\item The ordered products $1_{(1,l_1)}\cdots
  1_{(1,l_r)}1_{(0,d_1)\cdots 1_{(0,d_s)}}$ form a basis
  $\underline{\bf B}$ of $\alg_\P1$, and $H$ is injective.
\item The elements $1_Z$, for $Z \in \text{Irr}(\Higgs)$ constructed
  in section \ref{semicanbasis} form a basis
  $\bf B$ of $\alg_\P1$.
\end{enumerate}
\end{lem}
\begin{dem}
The first point follows easily from the definition of the convolution
product in each algebra.\\
To see (2), note that the generators of $\mathcal{H}_\P1$ are mapped into a
set of generators of $\bf H_\P1$ (see \cite{BK}).\\
For (3), we already know (see \ref{ord}) that these elements generates
$\mathcal{H}_\P1$ as a vector space. But these elements are mapped into
linearly independant elements of $\bf H_\P1$, and consequently form a
basis of $\mathcal{H}_\P1$. The injectivity of $H$ follows.

The point (4) follows from the construction of the elements in $\bf
B$. The (infinite) matrix which expresses these elements in the basis
$\underline{\bf B}$ is upper triangular with non-zero coefficient along
the diagonal.
\end{dem}

\section{Crystal structure}

The aim of this section is to endow the semicanonical basis with the
structure of a ``crystal'', in analogy with the case of quivers (see
\cite{L1},\cite{KS}). In the situation of quivers, the canonical basis
and the semi-canonical basis though often different  have isomorphic
crystals. In the situation of curves, the canonical basis (which is
constructed in \cite{Sc2}) is not known to have a crystal structure,
so the existence of such a structure for the semi-canonical basis is
interesting. 

\subsection{}
For a coherent sheaf $\CalF$ on $\P1$ and a line bundle
$\CalO(k)$, define 
\[
\text{rk}_k(\CalF)=\text{max}\{i\in \Nn| \ \text{Inj}(\CalO(k)^{\oplus i},\CalF)
\neq \emptyset\}
\]
where $\text{Inj}(\CalO(k)^{\oplus i},\CalF)$ is the set of injections in
$\Hom(\CalO(k)^{\oplus i},\CalF)$.\\
Now for $k \in \Zz$ and $s,n \in \Nn$, we define a locally closed
substack $\underline{\Lambda}^{(r,d)}_{k,n,s}$ of
$\underline{\Lambda}_\P1^{(r,d)}$ by
\[
\underline{\Lambda}^{(r,d)}_{k,n,s}=\{(\CalF,f) \in
\underline{\Lambda}^{(r,d)}_\P1, \  \text{rk}_k(\Ker(f))=s,\
\dim(\Hom(\CalO(k),\Ker(f)))=n \}
\]
In the following write $\alpha=(r,d)$, $\gamma=(s,sk)$ and $\beta=\alpha-gamma$, three elements of $K(\text{Coh}_\P1)$.\\
Define the following stack:
\[
\underline{\mathcal{E}}^{\alpha}_{k,n,s}=\{ (\CalF,f,i), \ (\CalF,f) \in
\underline{\Lambda}^{\alpha}_{k,n,s}, \  i \in
\text{Inj}(\CalO(k)^{\oplus s},\Ker(f)) \}
\]
representing the functor from the category of affine schemes to the
category of groupoids:
\begin{align*}
\Sigma\mapsto\{ (\CalF,f,i), \ \CalF
  \text{ is a coherent $\Sigma$-flat sheaf on $\P1\times \Sigma$, } \\
 & \hspace{-7cm} f \in \Hom(\CalF, \CalF \otimes \CalO_{\Sigma\times\P1}(-2)),
  i:\CalO_{\Sigma\times \P1}^s \rightarrow \Ker f, & \\
& \hspace{-8cm} \CalF_\sigma \text{ is of class $\alpha$ and $\text{rk}_k(\Ker f_\sigma)=s$},\ i_\sigma: \CalO(k)^s\hookrightarrow  \Ker f_\sigma \\
 & \hspace{-5cm} \text{for all closed point $\sigma \in \Sigma$}\}
\end{align*}
where a morphism $\psi$ between objects $(\CalF,f,i)$ and
$(\CalF',f',i')$ is an isomorphism $\psi:\CalF\simeq \CalF'$ such that
the following diagrams commute:
\[
\xymatrix{
 \CalF \ar[r]^\psi \ar[d]^f &  \CalF'\ar[d]^{f'} & \CalO_\Sigma(k)^s
 \ar[r]^i \ar[rd]^{i'} & \Ker f \ar[d]^\psi \\
\CalF(-2) \ar[r]^\psi & \CalF'(-2) &  & \Ker f'
}
\]
Now consider the following diagram:
\begin{equation}\label{diag}
\xymatrix{
 & \underline{\mathcal{E}}^{\alpha}_{k,n,s}  \ar[ld]_{p_1} \ar[rd]^{p_2}& \\
\underline{\Lambda}^{\alpha}_{k,n,s} & &
 \underline{\Lambda}^{\beta}
}
\end{equation}
where the maps $p_1$ and $p_2$ are defined on objects as follows:
\[
\begin{array}{cccc}
p_1: & (\CalF,f,i) & \mapsto & (\CalF,f) \\
p_2: & (\CalF,f,i)& \mapsto & (\CalF/i(\CalO(k)^s),f|_{\CalF/i(\CalO(k)^s)}),
\end{array}
\]
We begin with the following lemma:
\begin{lem}
$p_2(\underline{\mathcal{E}}^{\alpha}_{k,n,s})\subseteq \underline{\Lambda}_{k,n-s,0}^{\beta}$
\end{lem}
\begin{dem}
Take an element $(\CalF,i,f)$ in
$\underline{\mathcal{E}}^{\alpha}_{k,n,s}$ and define
$(\CalG,g)=p_2(\CalF,i,f)$. We have a short exact sequence:
\begin{equation}
\xymatrix{
0 \ar[r] & \CalO(k)^s \ar[r]^i & \CalF \ar[r]^j & \CalG \ar[r] & 0
}
\end{equation}
As $\Ext(\CalO(k),\CalO(k))=0$, we also have
\begin{equation}
\xymatrix{
0 \ar[r] & \Hom(\CalO(k),\CalO(k)^s) \ar[r]^{i' }& \Hom(\CalO(k),\CalF) \ar[r]^{j'} & \Hom(\CalO(k),\CalG) \ar[r] & 0
}
\end{equation}
Now the following diagram commutes for any $a\in\Hom(\CalO(k),\Ker f
)$
\begin{equation}
\xymatrix{
\CalO(k) \ar[r] \ar[rd]^a & \CalG \ar[r]^g & \CalG(-2) \\
 & \CalF \ar[u]^j \ar[r]^f & \CalF(-2) \ar[u]^j,
}
\end{equation}
so we see that $f\circ a=0$ implies $g\circ j'(a)=0$, then
$j'(\Hom(\CalO(k),\Ker f ))\subseteq \Hom(\CalO(k),\Ker g)$. In the
other way, if $g\circ j'(a)=0$, then as the kernel of the morphism
$\CalF(-2) \rightarrow \CalG(-2)$ is $\CalO(k-2)^s$, we have that
$\Im(f\circ a)\subseteq \CalO(k-2)^s$. The morphism $f\circ a$ lies
inside $\Hom(\CalO(k),\CalO(k-2)^s)=0$, so that $a\in
\Hom(\CalO(k),\Ker f)$. We just proved that we have a short exact sequence:
\[
0 \rightarrow \Hom(\CalO(k),\CalO(k)^s) \rightarrow \Hom(\CalO(k),\Ker
f) \rightarrow \Hom(\CalO(k),\Ker g) \rightarrow 0
\]
It follows that $\dim \Hom (\CalO(k),\Ker g)=n-s$.\\
Now we prove that there are no injections from $\CalO(k)$ into
$\Ker g$. Assuming that such an injection $h$ exists, consider an element
$h'\in \Hom(\CalO(k),\Ker f)$ such that $j'(h')=h$. Define
$h''=h'\oplus i \in \Hom(\CalO(k)\oplus \CalO(k)^s,\Ker f)$. From the
commutative diagram: 
\begin{equation*}
\xymatrix{
0 \ar[r] & \CalO(k)^s \ar[r]  & \CalO(k)^{s+1}
\ar[r]^{pr}\ar[d]^{h''}&  \CalO(k)\ar[d]^h \ar[r]& 0\\
 & & \CalF \ar[r]^j & \CalG & ,
}
\end{equation*}
we deduce that $pr(\Ker h'')=0$, i.e. $\Ker h''\subseteq \CalO(k)^s$. But the restriction of
$h''$ to $\CalO(k)^s$ is $i$, which is injective. We have proved that
$h''$ is an injection from $\CalO(k)^{s+1}$ into $\Ker f$, which is
impossible since $\text{rk}_k(\Ker f)=s$.
\end{dem}
The diagram \ref{diag} is then refined to the diagram:
\begin{equation}\label{diag2}
\xymatrix{
 & \underline{\mathcal{E}}^{\alpha}_{k,n,s}  \ar[ld]_{p_1} \ar[rd]^{p_2}& \\
\underline{\Lambda}^{\alpha}_{k,n,s} & &
 \underline{\Lambda}^{\beta}_{k,n-s,0}
}
\end{equation}
The main theorem of this section is the following:
\begin{thm}\label{CorrespIrr}
We have a natural bijection between irreducible components $Z_1$ of
$\underline{\Lambda}^{\alpha}_{k,n,s}$ and irreducible components $Z_2$
of $\underline{\Lambda}^{\beta}_{k,n-s,0}$. Under this
correspondence we have $\dim Z_1 = \dim Z_2 -2s(r-s)-s^2$.
\end{thm}
We have to study the maps $p_1$ and $p_2$. We start with $p_2$.
\begin{lem}
The map $p_2$ is smooth with connected fibers of dimension
$s(n-s)-2s(r-s)-s^2$.
\end{lem}
\begin{cor} The map $p_2$ induces a natural bijection between
  irreducible components $Z_2$ of
  $\underline{\Lambda}_{k,n-s,0}^{\beta}$ and irreducible components
  $Z_3$ of $\underline{\mathcal{E}}_{k,n,s}^{\alpha}$. Moreover, under
  this correspondence we have $\dim Z_3=\dim Z_2 -2s(r-s)+s(n-s)-s^2$.
\end{cor}
\begin{dem}
As usual, to prove the lemma we need to study locally the diagram
(\ref{diag}). To do so, define the locally closed subvariety
$S_{m,(k,n,s)}^{\alpha}$ of $S_m^{\alpha}$ by 
\[
S_{m,(k,n,s)}^{\alpha}=\{ (\phi,\CalF,f) \in S_m^{\alpha}\ | \
(\CalF,f)\in \underline{\Lambda}_{k,n,s}^{\alpha}\}
\]
For an integer $m<k$ set $d_1=\langle \CalO(m),\CalO(k)^{\oplus s}
\rangle=\dim \Hom(\CalO(m), \CalO(k)^{\oplus s})$, $d'= \langle
\CalO(m),\CalF\rangle= \dim \Hom(\CalO(m),\CalF)$ and
$d_2=d'-d_1$. Define
\begin{multline*}
E^{\alpha,\geq m}_{k,s,n}=\{(\phi,\CalF,f,i,h_1,h_2), \ (\phi,\CalF,f) \in
S_{m,(k,n,s)}^{\alpha},\\
i:\CalO(k)^s \hookrightarrow \Ker f ,\ h_1:\Cc^{d_1} \simeq
\Hom(\CalO(m),\CalO(k)^s),\\
h_2:\Cc^{d_2}\simeq \Hom(\CalO(m), \CalF/i(\CalO(k)^s)
\}
\end{multline*}
The group $G=GL_{d'}\times GL_{d_1} \times GL_{d_2}$ acts
naturally on $E^{\alpha,\geq m}_{k,s,n}$ and the quotient stack is
$\underline{\mathcal{E}}_{k,n,s}^{\alpha,\geq m}$.\\
Introduce $C=\{(V,a,b)| V \subseteq \Cc^{d'},  a:V \simeq \Cc^{d_1}, b: \Cc^{d'}/V\simeq \Cc^{d_2}\}$.\\
We define $q_2$ as follows:
\[
\begin{array}{cccc}
q_2: & E^{\alpha,\geq m}_{k,n,s} & \rightarrow &
 S_{m,(k,n-s,0)}^{\beta}\times S_{m,(k,s,s)}^{\gamma} \times
 C \\
 & (\phi,\CalF,f,i,h_1,h_2) & \mapsto &
 ((\psi_1,\CalG,g),(\psi_2,\CalO(k)^{\oplus s},0),(V,a,b))
\end{array}
\]
where:
\begin{enumerate}
\item $\CalG:=\CalF/i(\CalO(k)^s)$,
\item $\psi_1: \CalO(m)^{\oplus d_2} \twoheadrightarrow \CalG$ is deduced from
  $\phi$ and $h_2$,
\item $\psi_2:\CalO(m)^{\oplus d_1} \twoheadrightarrow
  \CalO(k)^{\oplus s}$ is deduced from $\phi$ and $h_1$,
\item $(V,a,b)\in Z$ is defined by $V=\phi_*(\Hom
  (\CalO(m),\CalO(k)^s))\subseteq
  \Cc^{d_1+d_2}=\phi_*(\Hom(\CalO(m),\CalF))$ (via $i$) and $a$ and
  $b$ are deduced from the diagram
\begin{equation*}
\xymatrix{
0 \ar[r] & \Hom(\CalO(m),\CalO(k)^s) \ar[r]^i\ar[d]_{h_1}& \Hom(\CalO(m),\CalF)
\ar[r]\ar[d]_{\phi_*}& \Hom(\CalO(m),\CalG) \ar[r]\ar[d]_{h_2} & 0 \\
0 \ar[r]& \Cc^{d_1} \ar[r] &\Cc^{d'} \ar[r] &\Cc^{d_2}\ar[r]& 0
}
\end{equation*}
\end{enumerate}
\begin{lem}
The map $q_2$ is an affine fibration, with fibers of dimension $-\langle
\gamma,\beta\rangle -\langle \beta,\gamma \rangle+ d_1d_2 +s(n-s)$.
\end{lem}
\begin{dem}
Let us fix some notations. We will write $\mathcal{L}=\CalO(k)^s$, $\CalO_\mathcal{L}=\CalO(m)^{ d_1}$,
  $\CalO_\CalG=\CalO(m)^{ d_2}$, $\CalO_\CalF=\CalO(m)^{ d'}$.\\
For maps
\begin{equation*}\begin{cases}
\psi_1:\CalO(m)^{d_2}\twoheadrightarrow \CalG \\
\psi_2:\CalO(m)^{d_1}\twoheadrightarrow \CalO(k)^s \\
\phi:\CalO(m)^{d'}\twoheadrightarrow \CalF 
\end{cases}
\end{equation*}
write $\mathcal{K}_\mathcal{L}$, $\mathcal{K}_\CalG$ and
  $\mathcal{K}_\CalF$ for the corresponding kernels. We denote
  $i_\mathcal{L} :\mathcal{K}_\mathcal{L} \hookrightarrow
  \CalO_\mathcal{L}$ and $i_\CalG :\mathcal{K}_\CalG \hookrightarrow
  \CalO_\CalG$ the corresponding injections.\\
Let us describe the morphism $q_2$ in terms of some diagrams. Elements
  in the space $E_{k,n,s}^{\alpha,\geq n}$ are in canonical bijection
  with commutative diagrams
\begin{equation}\label{elementfiber}
\xymatrix{
  0\ar[d] & 0\ar[d] & \\
 \CalO_\mathcal{L} \ar[d]^{a'}\ar[r]^\xi&\mathcal{L}\ar[r]\ar[d]^i&0 \\
 \CalO_\CalF\ar[d]^{b'}
 \ar[r]^\phi&\mathcal{F}\ar[r]\ar[d] & 0\\
 \CalO_\CalG
 \ar[d]\ar[r]^\eta&\mathcal{G}\ar[d]\ar[r]&0 \\
 0 & 0 & 
}
\end{equation}
together with a map $f\in \Hom(\CalF,\CalF(-2))$ such that
$i:\mathcal{L} \hookrightarrow \Ker f$. Indeed, the maps $\xi,\eta$
are deduced from $h_1,h_2$ by the formulas:
\begin{align*}
\xi=\text{can} \circ h_1,\\
\eta=\text{can} \circ h_2
\end{align*}
where can is the evaluation map and $a',b'$ are defined uniquely so as
to make \ref{elementfiber} commute. Recall that in the construction of
$\text{Hilb}_{\CalO(m)^{d'},\alpha}$, two maps
$\phi:\CalO(m)^{d'}\twoheadrightarrow \CalF$, $\phi':\CalO(m)^{d'}
\twoheadrightarrow \CalF'$ are equivalent if $\Ker \phi=\Ker
\phi'$. We use the same equivalence relation for diagrams.\\
Similarly, points in $S_{m,(k,n-s,0)}^{\beta}\times S_{m,(k,s,s)}^{\gamma} \times
 C$ correspond bijectively to diagrams
\begin{equation}\label{elementbottom}
\xymatrix{
  0 \ar[d]& & \\
 \CalO_\mathcal{L} \ar[d]^{a'}\ar[r]^{\psi_2}&\mathcal{L}\ar[r]&0 \\
 \CalO_\CalF\ar[d]^{b'}
 & & \\
 \CalO_\CalG
 \ar[d]\ar[r]^{\psi_1}&\mathcal{G}\ar[r]&0 \\
0 & &
}
\end{equation}
together with an element $g\in \Hom(\CalG,\CalG(-2))$.\\
The horizontal sequences are the elements of
$S_{m,(k,n-s,0)}^{\beta}$ and $S_{m,(k,s,s)}^{\gamma}$, and the
vertical sequence is deduced from $(V,a,b)$ by
\begin{equation*}
\xymatrix{
0 \ar[r] & V \ar[r] \ar[d]^a & \Cc^{d_1+d_2} \ar[r] \ar@2{{}-{}}[d]& \Cc^{d_1+d_2}/V \ar[d]^b \ar[r]& 0\\
0 \ar[r] & \Cc^{d_1} \ar[r] & \Cc^{d_1+d_2} \ar[r]& \Cc^{d_2}  \ar[r]& 0
}
\end{equation*}
and then tensoring by
$\CalO(m)$.\\
The map $q_2$ assigns to a diagram as in \ref{elementfiber} its
subdiagram \ref{elementbottom} and the element $g$ deduced from
$f$. We may complete the diagrams \ref{elementfiber},
\ref{elementbottom} by adding kernels of $\xi,\eta,\phi$:
\begin{equation}\label{elementfiber2}
\xymatrix{
  & 0\ar[d] & 0\ar[d] & 0\ar[d] & \\
 0\ar[r]& \mathcal{K}_\mathcal{L}\ar[r]^{i_\mathcal{L}}\ar[d] &
  \CalO_\mathcal{L} \ar[d]^{a'}\ar[r]^\xi &\mathcal{L}\ar[r]\ar[d]^i&0 \\
 0 \ar[r]& \mathcal{K}_\CalF \ar[r]^{i_\CalF} \ar[d]& \CalO_\CalF\ar[d]^{b'}
 \ar[r]^{\phi}&\mathcal{F}\ar[r]\ar[d] & 0\\
 0\ar[r] & \mathcal{K}_\CalG \ar[r]^{i_\CalG} \ar[d] & \CalO_\CalG
 \ar[d]\ar[r]^\eta&\mathcal{G}\ar[d]\ar[r]&0 \\
  & 0 &0 & 0 & 
}
\end{equation}
\begin{equation}\label{elementbottom2}
\xymatrix{
  &  & 0\ar[d] &  & \\
 0\ar[r]& \mathcal{K}_\mathcal{L}\ar[r]^{i_\mathcal{L}} &
  \CalO_\mathcal{L} \ar[d]^{a'}\ar[r]^\xi &\mathcal{L}\ar[r]&0 \\
  & & \CalO_\CalF\ar[d]^{b'}
 & & \\
  0\ar[r] & \mathcal{K}_\CalG \ar[r]^{i_\CalG}  & \CalO_\CalG
 \ar[d]\ar[r]^\eta&\mathcal{G}\ar[r]&0 \\
  &  &0 &  & 
}
\end{equation}
Let us fix a point $x=(\psi_1,\CalG,g,\psi_2,\mathcal{L},V,a,b)\in
S_{m,(k,n-s,0)}^{\beta}\times S_{m,(k,s,s)}^{\gamma} \times C$, and
denote the fiber $F=q_2^{-1}(x)$.\\
We will use the following lemma:
\begin{lem}\label{formkernel}
For any $m\in \Zz$ and any $\alpha \in K(\text{Coh}_\P1)$, the kernel
of any map $(\phi:\CalO(m)^{d(m,\alpha)} \twoheadrightarrow \CalF)\in
Q_m^\alpha$ is isomorphic to
$\CalO(m-1)^{d(m,\alpha)-\text{rk}(\alpha)}$.\\
Conversely, for any
embedding $\CalO(m-1)^{d(m,\alpha)-\text{rk}(\alpha)} \subseteq
\CalO(m)^{d(m,\alpha)}$, the map $\phi: \CalO(m)^{d(m,\alpha)}
\twoheadrightarrow
\CalO(m)^{d(m,\alpha)}/\CalO(m-1)^{d(m,\alpha)-\text{rk}(\alpha)}$
belongs to $Q_m^\alpha$.
\end{lem}
\begin{dem}
 The kernel $\Ker \phi$ is a
  vector bundle (as a subsheaf of a vector bundle) and we write it as
  $\oplus_{i=1}^{d(m,\alpha)-\text{rk}(\alpha)} \CalO(k_i)$.
The morphism 
\[
\phi_*:
\Cc^{d(m,\alpha)}=\Hom(\CalO(m),\CalO(m)^{d(m,\alpha)}) \rightarrow
  \Hom(\CalO(m), \CalF)
\]
is an isomorphism, so the kernel
$\Hom(\CalO(m),\Ker \phi)$ is zero. Then we have $k_i\leq m-1$ for
any $i$.\\
Now we know that the degree of $\Ker \phi$ is equal to
$md(m,\alpha)-\text{deg}(\alpha)=\sum_i k_i$. We have:
\[
\sum_i k_i=md(m,\alpha)-\text{deg}(\alpha)) \leq (m-1)(d(m,\alpha)-\text{rk}(\alpha))
\]
But as $d(m,\alpha)=(1-m)\text{rk}(\alpha)+\text{deg}(\alpha)$, we
have that the right-hand side is $md(m,\alpha)-(d(m,\alpha)+(m-1)\text{rk}(\alpha))=md(m,\alpha)-\text{deg}(\alpha)$
and this inequality is in fact an equality, so that each $k_i$ is
equal to $m-1$.\\
For the converse, if $\Ker \phi\simeq
\CalO(m-1)^{d(m,\alpha)-\text{rk}(\alpha)}$ then the morphism
\[
\Hom(\CalO(m),\CalO(m)^{d(m,\alpha)}) \xrightarrow{\phi_*}
\Hom(\CalO(m),\CalF)
\] is surjective because
$\Ext(\CalO(m),\CalO(m-1))=0$. As it is also injective, $\phi_*$ is an
isomorphism $\Cc^{d(m,\alpha)} \simeq \Hom(\CalO(m),\CalF)$.
\end{dem}
The set of classes of maps $\phi:\CalO_\CalF \twoheadrightarrow \CalF$
making \ref{elementfiber} commutative is in bijection with the set of
subsheaves $\mathcal{K}_\CalF \subseteq \CalO_\CalF$ satisfying:
\begin{equation}\label{condkernel}
\begin{cases}
\mathcal{K}_\CalF \cap
a'(\CalO_\mathcal{L})=a'(\mathcal{K}_\mathcal{L})\\
b'(\mathcal{K}_\CalF)=\mathcal{K}_\CalG
\end{cases}
\end{equation}
Indeed, if \ref{condkernel} holds then $\mathcal{K}_\CalF$ fits
in a short exact sequence
\[
0\rightarrow \mathcal{K}_\mathcal{L}\rightarrow \mathcal{K}_\CalF
\rightarrow \mathcal{K}_\CalG \rightarrow 0
\]
hence $\mathcal{K}_\CalF\simeq
\CalO(m-1)^{d(m,\alpha)-\text{rk}(\alpha)}$ and we apply the second
part of lemma \ref{formkernel}.\\
Subsheaves $\mathcal{K}_\CalF \subseteq \CalO_\CalF$ satisfying
\ref{condkernel} form a principal
$\Hom(\mathcal{K}_\CalG,\mathcal{L})$-space. Indeed
\begin{multline*}
\{\mathcal{K}_\CalF \subseteq \CalO_\CalF\ |\ \text{\ref{condkernel} is
  satisfied}\}= \{ \mathcal{K}_\CalF' \subseteq
  \CalO_\CalF/\mathcal{K}_\mathcal{L}\ | \ \mathcal{K}_\CalF'\cap
  \mathcal{L}=0,b'(\mathcal{K}_\CalF')=\mathcal{K}_\CalG\} \\
=\{s: \mathcal{K}_\CalG \rightarrow
  \CalO_\CalF/\mathcal{K}_\mathcal{L}\ |\  b'\circ s=Id_{\mathcal{K}_\CalG}\}
\end{multline*}
and if $s,s'$ are two sections $\mathcal{K}_\CalG \rightarrow
\CalO_\CalF/ \mathcal{K}_\mathcal{L}$ as above then $s-s' \in
\Hom(\mathcal{K}_\CalG,\CalO_\mathcal{L}/\mathcal{K}_\mathcal{L})=\Hom(\mathcal{K}_\CalG,\mathcal{L})$.\\
For convenience, let us chose a section $s_0$ as abive. This
corresponds to an identification
$\CalO_\CalF/\mathcal{K}_\mathcal{L}\simeq \mathcal{L} \oplus
\CalO_\CalG$. Then to $u\in \Hom(\mathcal{K}_\CalG,\mathcal{L})$ we
associate the diagram
\begin{equation}
\xymatrix{
  & 0\ar[d] & 0\ar[d] & 0\ar[d] & \\
 0\ar[r]& \mathcal{K}_\mathcal{L}\ar[r]^{i_\mathcal{L}}\ar[d] &
  \CalO_\mathcal{L} \ar[d]^{a'}\ar[r]^{\psi_2} &\mathcal{L}\ar[r]\ar[d]^i&0 \\
 0 \ar[r]& \mathcal{K}_\CalF \ar[r]^{i_\CalF} \ar[d]& \CalO_\CalF\ar[d]^{b'}
 \ar[r]^{\phi_u}&\mathcal{F}\ar[r]\ar[d] & 0\\
 0\ar[r] & \mathcal{K}_\CalG \ar[r]^{i_\CalG} \ar[d] & \CalO_\CalG
 \ar[d]\ar[r]^{\psi_1}&\mathcal{G}\ar[d]\ar[r]&0 \\
  & 0 &0 & 0 & 
}
\end{equation}
Note that $\CalF\simeq\text{Coim}(\mathcal{K}_\CalG
\hookrightarrow^{(i_\CalG,u)} \CalO_\CalG \oplus \mathcal{L})$.\\
 It remains to describe the possible choices for the map $f$ in the fiber. Such an element
 verifies two conditions:

(*) $f|_\mathcal{L}=O$

(**) $f'|_\CalG=g$, where $f' \in \Hom(\CalG,\CalG(-2))$ is deduced
from $f$.\\
We have the short exact sequence derived from $u$:
\[
0 \rightarrow \mathcal{L} \rightarrow \CalF \rightarrow \CalG
\rightarrow 0
\]
From the following commutative diagram, where vertical maps are
deduced from Serre duality
\begin{equation*}
\xymatrix{
 0 \ar[r]& \Ext(\CalF,\CalG)^* \ar[r] \ar@2{{}-{}}[d] & \Ext(\CalF,\CalF)^* \ar@2{{}-{}}[d]\ar[r]& \Ext(\CalF,\mathcal{L})^*\ar@2{{}-{}}[d] \\
0 \ar[r]& \Hom(\CalG,\CalF(-2)) \ar[r] & \Hom(\CalF,\CalF(-2)) \ar[r]^{|_\mathcal{L}} & \Hom(\mathcal{L},\CalF(-2))
}
\end{equation*}
We see that $f_{|\mathcal{L}}=0$ is equivalent to $f\in \Ext(\CalF,\CalG)^*$. The other condition is given by:
\begin{equation*}
\begin{array}{ccccccccc}
0 & \rightarrow & \Ext(\mathcal{L},\CalG)^* & \rightarrow & \Ext(\CalF,\CalG)^* & \rightarrow & \Ext(\CalG,\CalG)^* & \xrightarrow{\theta_u} & \Hom(\mathcal{L},\CalG)^*\\
 & & & & f & \mapsto & g & &
\end{array}
\end{equation*}
where $\theta_u$ is the connecting morphism.\\
The possible choices of $f$ in the fiber is then a principal $\Ext(\mathcal{L},\CalG)^*$-space, when we have the condition $\theta_u(g)=0$.\\
To sum up, we have shown that the fiber $F$ is isomorphic to the
subspace of pairs $(u,v)\in \Hom(\mathcal{K}_\CalG,\mathcal{L})\oplus
\Ext(\mathcal{L},\CalG)^*$ satisfying $\theta_u(g)=0$. We need to
describe more precisely the map $\theta_u$.
\begin{lem}
The map $\theta_u:\Ext(\CalG,\CalG)^* \rightarrow
\Hom(\mathcal{L},\CalG)^*$ is given by
\[
\theta_u(g)(h)= a_g(h\circ u)
\]
for any $h\in \Hom(\mathcal{L},\CalG)$, where $a_g$ is the image
of $g$ in $\Hom(\mathcal{K}_\CalG,\CalG)^*$.
\end{lem}
\begin{dem}
We claim that the following diagram is commutative
\begin{equation*}
\xymatrix{
\Ext(\CalG,\CalG)^* \ar[r]^{\theta_u} \ar@{^{(}->}[d] & \Hom(\mathcal{L},\CalG)^* \ar@{^{(}->}[d] \\
\Hom(\mathcal{K}_\CalG,\CalG)^* \ar[r]^{\theta_u'} & \Hom(\CalO_\CalG\oplus \mathcal{L},\CalG)^*
}
\end{equation*}
where $\theta_u'$ is induced by the injection $\mathcal{K}_\CalG
\xrightarrow{(i_\mathcal{G},u)} \CalO_\CalG \oplus \mathcal{L}$. \\
To see this, apply $\Hom(.,\CalG)$ to the diagram
\[
\xymatrix{ & & \mathcal{L}\ar[d]\ar[dr] & & \\
0 \ar[r]&\mathcal{K}_\CalG \ar[r]\ar@{=}[d] &\CalO_\CalG \oplus \mathcal{L} \ar[r] \ar[d]
& \CalF \ar[r] \ar@{->>}[d] & 0\\
0\ar[r] & \mathcal{K}_\CalG \ar[r] & \CalO_\CalG \ar[r] & \CalG \ar[r]& 0
}
\]
to get the construction of the connecting morphism $\theta^*_u$
\[
\xymatrix{
0 \ar[r] & \Hom(\CalG,\CalG) \ar[r]\ar[d]& \Hom(\CalF,\CalG) \ar[r]\ar[d]
  & \Hom(\mathcal{L},\CalG) \ar@{=}[d] \ar@{-->}[dl]\ar `[rr]`[rrddd][dddll]^{\theta_u^*}&  &\\
0 \ar[r] & \Hom(\CalO_\CalG,\CalG) \ar[r]\ar[d] & \Hom(\CalO_\CalG \oplus
\mathcal{L},\CalG) \ar@{-->}[dl]\ar@<1ex>[r]^-{proj} \ar[d]^{(i_\mathcal{G},u)} &
\Hom(\mathcal{L},\CalG) \ar[r]\ar@<1ex>[l]^-{can} & 0 &\\
 & \Hom(\mathcal{K}_\CalG,\CalG)\ar[d] \ar@{=}[r] &
 \Hom(\mathcal{K}_\CalG,\CalG) & && \\
& \Ext(\CalG,\CalG) & & & &
}
\]
as the composition of the two dotted arrows and the map
$\Hom(\mathcal{K}_\CalG,\CalG) \rightarrow \Ext(\CalG,\CalG)$, which is exactly the dual
of our claim.\\
So for $h\in \Hom(\mathcal{L},\CalG)$, we have $\theta_u'(a_g)(h)=a_g(h \circ (i_\CalG,u))$. Then $\theta_u$ is obtained by evaluating $\theta_u$ on the projection of $h \circ (i_\CalG, u)$ into $\Hom(\mathcal{L},\CalG)$, i.e. $\theta_u(g)(h)=a_g(h\circ u)$.
\end{dem}
We can now consider a new linear map $\theta$ defined from $\theta_u$, this time considering the dependence on $u$:
\[
\begin{array}{cccc}
\theta: & \Hom(\mathcal{K}_\CalG,\mathcal{L})  & \rightarrow & \Hom(\mathcal{L},\CalG)^*\\
 & u & \mapsto & \theta_u(g)
\end{array}
\]
We have proved the following statement:
\begin{lem}\label{fiber}
The fiber $F$ is isomorphic to $\Ext(\mathcal{L},\CalG)^* \oplus \Ker \theta$.
\end{lem}
It remains to give the dimension of the fiber $F$. We need the following lemma:
\begin{lem}\label{dim}
We have $(\Im \theta)^\bot= \Hom(\mathcal{L},\Ker g)$.
\end{lem}
\begin{dem}
Take $h\in \Hom(\mathcal{L},\CalG)$, and define $\CalI:=\Im h$. Now
$h\in (\Im\theta)^\bot$ is equivalent to:
\[
\forall u\in \Hom(\mathcal{K}_\CalG,\mathcal{L}),\ a_g(h\circ u)=0
\]
We have a natural map $\Hom(\mathcal{K}_\CalG,\CalI) \xrightarrow{p}
\Ext(\CalG,\CalI)$, and by definition of $a_g$, we have
$a_g(v)=g(p(v))$ for any $v\in \Hom(\mathcal{K}_\CalG,\CalI)$.
\begin{equation*}
\xymatrix{
\Hom(\mathcal{K}_\CalG,\CalI) \ar[r]^p \ar[d]^{a_g} &
\Ext(\CalG,\CalI) \ar[d]^g \\
\Cc & \Cc
}
\end{equation*}
From the surjection $h: \mathcal{L} \twoheadrightarrow \CalI$, we
have a surjection $\Ext(\CalG,\mathcal{L})\twoheadrightarrow
\Ext(\CalG,\CalI)$, and we can make use of the following commutative
diagram:
\begin{equation*}
\xymatrix{
\Hom(\mathcal{K}_\CalG,\mathcal{L}) \ar[r]^{h_*'} \ar[d]^{p'} &
\Hom(\mathcal{K}_\CalG,\CalI) \ar[d]^p &
\\
\Ext(\CalG,\mathcal{L}) \ar[d] \ar[r]^{h_*}& \Ext(\CalG,\CalI) \ar[r] & 0\\
 0 & &
}
\end{equation*}
So that we have the following chain of equivalence:
\[
\begin{array}{ccc}
h \in (\Im \theta)^\bot & \Leftrightarrow & \forall u\in
\Hom(\mathcal{K}_\CalG,\mathcal{L}),\ a_g(h\circ u)=0=g(p(h \circ u)) \\
& \Leftrightarrow & g(p(h_*'(\Hom(\mathcal{K}_\CalG, \mathcal{L}))))=0\\
& \Leftrightarrow & g(h_*(p'(\Hom(\mathcal{K}_\CalG,\mathcal{L}))))=0 \\
& \Leftrightarrow & g(h_*(\Ext(\CalG,\mathcal{L})))=0  \text{  (by
  surjectivity of }p')\\
& \Leftrightarrow & g|_{\Ext(\CalG,\CalI)}=0 \text{  (by surjectivity of }h_*)
\end{array}
\]
Now the restriction $g|_{\Ext(\CalG,\CalI)}$ is equal to the image of
$g$ by the morphism $\Ext(\CalG, \CalG)^* \rightarrow
\Ext(\CalG,\CalI)^*$ which is just the restriction morphism, as we see
from Serre duality:
\[
\xymatrix{
 \Ext(\CalG,\CalG)^* \ar[r] \ar@2{{}-{}}[d] & \Ext(\CalG,\CalI)^*
 \ar@2{{}-{}}[d]\\
\Hom(\CalG,\CalG(-2)) \ar[r]^{|_\CalI} & \Hom(\CalI,\CalG(-2))
}
\]
so that $g|_{\Ext(\CalG,\CalI)}=g|_\CalI$, this time considered as an
element of $\Hom(\CalG,\CalG(-2))$.\\
We have proved that $h\in(\Im \theta)^\bot \Leftrightarrow g|_\CalI=0$, which
by definition is equivalent to $\CalI \subseteq \Ker g$.
\end{dem}
Lemma \ref{fiber} gives that $q_2$ is an affine fibration with
connected fibers. Lemma
\ref{dim} allows us to compute the dimension of the fiber. As $\dim
\Hom (\mathcal{L},\Ker g)=s \dim \Hom(\CalO(k),\Ker g)=s(n-s)$, we have:
\begin{align*}
\dim \Ker \theta &= &\dim
\Ext(\mathcal{L},\CalG) + \dim
\Hom(\mathcal{K}_\CalG,\mathcal{L})-(\dim \Hom(\mathcal{L},\CalG)
-s(n-s))\\
&=&\dim \Ext(\mathcal{L},\CalG)- \dim
\Hom(\mathcal{L},\CalG)+\dim\Hom(\mathcal{K}_\CalG,\mathcal{L}) +s(n-s)\\ 
&=&-\langle\mathcal{L},\CalG\rangle +\langle
\mathcal{K}_\mathcal{L}, \CalG\rangle +s(n-s)\\
&=&-\langle \mathcal{L},\CalG\rangle +\langle \CalO_\CalG,\mathcal{L}
\rangle-\langle \CalG,\mathcal{L} \rangle +s(n-s)\\
&=&-\langle \gamma,\beta\rangle -\langle \beta,\gamma
\rangle +d_1d_2 +s(n-s)
\end{align*}
In the above calculation we have used
$\Ext(\mathcal{K}_\CalG,\mathcal{L})=0$ and
$\langle \CalO_\CalG,\mathcal{L}\rangle=\langle
\CalO(m)^{d_2},\mathcal{L}\rangle=d_2 \langle \CalO(m),\mathcal{L}
\rangle=d_2d_1$.
\end{dem}
As the map $q_2$ is $G$-equivariant, we can pass to the quotient
to obtain a map $q_2'$, which is also an affine fibration with
connected fibers:
\[
q_2':\underline{\mathcal{E}}_{k,n,s}^{\alpha,\geq m} \rightarrow
\underline{\Lambda}_{k,n-s,0}^{\beta,\geq m} \times
\underline{\Lambda}_{k,s,s}^{\gamma,\geq m}\times [C/GL_{d'}]
\]
The variety $C$ is an homogenous $GL_{d'}$-variety (hence smooth) of
dimension $d'^2-d_1d_2$, so the quotient is a smooth (connected) stack of
dimension $-d_1d_2$.\\
We have the following diagram:
\[
\xymatrix{
\underline{\mathcal{E}}_{k,n,s}^{\alpha,\geq m} \ar[rr]^-{q_2'}
\ar[rrd]^{p_2} & & \underline{\Lambda}_{k,n-s,0}^{\beta,\geq m}\times
\underline{\Lambda}_{k,s,s}^{\gamma,\geq m}\times [C/GL_{d'}] \ar[d]^{proj}\\
 &  & \underline{\Lambda}_{k,n-s,0}^{\beta,\geq m}
}
\]
As the stack $\underline{\Lambda}_{k,s,s}^{(s,sk)}$ is smooth
connected of dimension
$-\langle \gamma,\gamma \rangle=-s^2$, the morphism $p_2$ is smooth
with connected fibers of dimension $-\langle \beta, \gamma \rangle
-\langle \gamma, \beta \rangle - \langle \gamma,\gamma \rangle+s(n-s) =-2s(r-s)+s(n-s)-s^2$.
\end{dem}
Now we study the map $p_1$.
\begin{prop}
There is a natural bijection between between irreducible components $Z_1$ of
$\underline{\Lambda}^{\alpha}_{k,s,n}$ and irreducible components $Z_3$
of $\underline{\mathcal{E}}_{k,n,s}^{\alpha}$. Under this
correspondence we have $\dim Z_3=\dim Z_1 +s(n-s)$.
\end{prop}
\begin{dem}
We enlarge our stack
$\underline{\mathcal{E}}_{k,n,s}^{\alpha}$. Define a stack classifying
isomorphism classes of objects:
\[
\underline{F}_{k,n,s}^{\alpha}=\{(\CalF,f,i) |\ 
(\CalF,f)\in\underline{\Lambda}_{k,n,s}^{\alpha} , \ i \in
Gr_s^{\Hom(\CalO(k),\Ker f)} \}
\]
where as usual morphisms between objects $(\CalF,f,i)$ and
$(\CalF',f',i')$ are isomorphisms $\psi:\CalF \simeq \CalF'$ such that
the following diagrams commute:
\[
\xymatrix{
 \CalF \ar[r]^\psi \ar[d]^f &  \CalF'\ar[d]^{f'} & \CalO(k)^s
 \ar[r]^i \ar[rd]^{i'} & \Ker f \ar[d]^\psi \\
\CalF(-2) \ar[r]^\psi & \CalF'(-2) &  & \Ker f'
}
\]
The substack $\underline{\mathcal{E}}_{k,n,s}^{\alpha}$ is easily seen
to be an open dense substack of $\underline{F}_{k,n,s}^{(r,d)}$ (as
the condition $i$ injective is open
in the irreducible variety $Gr_s^{\Hom(\CalO(k),\Ker f)}$), and the map
  $p_1$ naturally extends to $\underline{F}_{k,n,s}^{\alpha,\beta}$.\\
Define the stack:
\[\underline{G}_{k,n,s}^{(r,d)}=\{ (\CalF,f,i,h)| \ (\CalF,f,i)
\in \underline{F}_{k,n,s}^{(r,d)}, \ h:\Cc^n\simeq \Hom(\CalO(k),\Ker
f) \}
\]
and
\[\underline{G'}_{k,n,s}^{(r,d)}=\{ (\CalF,f,i,h)| \ i \in
Gr_s^n\}
\]with the rest of the data as in
$\underline{G}_{k,n,s}^{(r,d)}$. We have natural maps which lead to the
following commutative diagram:
\[
\xymatrix{
\underline{G}_{k,n,s}^{(r,d)} \ar[r]^{\nu_1} \ar[d]_t &
\underline{F}_{k,n,s}^{(r,d)} \ar[r]^{p_1} &
\underline{\Lambda}_{k,n,s}^{(r,d)} \\
\underline{G'}_{k,n,s}^{(r,d)} \ar[r]^{\nu_2} &
\underline{\Lambda}_{k,n,s}^{(r,d)}\times Gr_n^s \ar[ur]_{proj} &
}
\]
where we have:
\begin{enumerate}
\item The map $t$ is defined by $t(\CalF,f,i,h)=(\CalF,f,i',h)$ where $i' \in Gr_s^n$ is
  deduced from $i$ via $h$. It is clearly an isomorphism.
\item The maps $\nu_1$ and $\nu_2$ (defined by
  $\nu(\CalF,f,i,h)=(\CalF,f,i)$) are $GL_n$ principal bundles.
\end{enumerate}
Consequently, this diagram induces a bijection between
Irr($\underline{F}_{k,s,n}^{(r,d)}$) and Irr($\underline{\Lambda}_{k,s,n}^{(r,d)}$). But as
$\underline{\mathcal{E}}_{k,s,n}^{(r,d)}$ is an open dense substack of
$\underline{F}_{k,s,n}^{(r,d)}$, it
also gives a bijection between
Irr($\underline{\mathcal{E}}_{k,s,n}^{(r,d)}$) and
Irr($\underline{\Lambda}_{k,s,n}^{(r,d)}$). Moreover, under this
correspondence $Z_1 \leftrightarrow Z_2$ we have $\dim Z_1 =\dim Z_2
+s(n-s)$.
\end{dem}
Now Theorem \ref{CorrespIrr} is a consequence of the two preceding
lemmas.
\subsection{}
As a corollary of theorem \ref{CorrespIrr}, we state
\begin{thm}\label{dim2}
We have
\begin{enumerate}
\item For $k \in \Zz$, $n,s\in \Nn$, the stack
  $\underline{\Lambda}_{k,n,s}^{(r,d)}$ is either empty or pure of
  dimension $\dim T^*\Coh_\P1^{(r,d)}/2=-r^2$
\item $\underline{\Lambda}_\P1^{(r,d)}$ is a closed substack of 
  $T^*\Coh_\P1^{(r,d)}$ of pure dimension $\dim T^*\Coh_\P1^{(r,d)}/2=-r^2$.
\end{enumerate}
\end{thm}
\begin{dem}
The proof is based on the following easy lemma:
\begin{lem}\label{rec}If $r>0$, we have
\[
\underline{\Lambda}_\P1^{(r,d)}= \bigcup_{k,n,s>0}
\underline{\Lambda}_{k,n,s}^{(r,d)},
\]
and the sum on the right-hand side is locally finite.
\end{lem}
\begin{dem}
It is a consequence of $\text{rk}(\Ker f) >0$ when
$\text{rk}(\CalF)>0$ (see also proof of lemma \ref{connected1}).
\end{dem}
Now we proceed by induction on $r$.\\
For $r=0$, the result follows from \cite{La1}, theorem 3.3.13.\\
For $r>0$, consider first the case $s>0$. \\
Then it follows from \ref{CorrespIrr} that irreducible
components of $\underline{\Lambda}_{k,n,s}^{(r,d)}$ are in
correspondance with irreducible components of
$\underline{\Lambda}_{k,n-s,0}^{(r-s,d-sk)}$, which from the induction
hypothesis are of dimension $-(r-s)^2$. It follows that
$\underline{\Lambda}_{k,n,s}^{(r,d)}$ is pure of dimension
$-(r-s)^2-2s(r-s)-s^2=-r^2$ for $s>0$.

By lemma \ref{rec}, $\underline{\Lambda}_\P1^{(r,d)\geq m}$ is a
finite union of such $\underline{\Lambda}_{k,n,s}^{(r,d)\geq m}$ with
$s>0$ for any $m\in \Zz$;
it follows that $\underline{\Lambda}_\P1^{(r,d)\geq m}$, and hence
$\underline{\Lambda}_\P1^{(r,d)}$ is pure of dimension $-r^2$.

To finish the induction step it remains to prove that $\underline{\Lambda}_{k,n,0}^{(r,d)}$ is pure
of the right dimension. Unfortunately this subspace is not open in
$\underline{\Lambda}_\P1^{(r,d)}$, so we have to work a little more.

We will use a refinement of lemma \ref{Dec}. Consider the substack
$\underline{L}^{r,d,l}_{k,n,0}$ of $\underline{L}^{r,d,l}$ defined as
(cf section \ref{Irredcomp})
\begin{multline*}
\underline{L}^{r,d,l}_{k,n,0}=\{(\mathcal{V},\tau,f_1,f_2,f_3)\in \underline{L}^{r,d,l}|
\dim\Hom(\CalO(k), \Ker f)=n, \\
\text{rk}_k(\Ker f)=0\}.
\end{multline*}
 Restrict the
diagram \ref{Dec} to the substack
$\underline{L}^{r,d,l}_{k,n,0}$ to obtain:
\[
\xymatrix{
 & \underline{L}^{r,d,l}_{k,n,0} \ar[dl]_{\pi_1} \ar[dr]^{\pi_2} & \\
\underline{\Lambda}^{r,d,l}_{k,n,0} & & T^*\Vec^{r,d-l}_{k,0,0}
 \times \underline{\Lambda}_{k,n,0}^{0,l}
}
\]
where $T^*\Vec^{r,d-l}_{k,0,0}=\{(\CalF,f)\in T^*\Vec^{r,d-l} |
\Hom(\CalO(k),\Ker f)=0\}$ is an open substack of
$T^*\Vec_\P1^{r,d-l}$. Indeed, for a vector bundle $\CalF=\oplus_i
\CalO(k_i)$, $k_i \geq k_{i+1}$, we always have $\CalO(k_1)\subseteq
\Ker f $ so that the condition $\Hom(\CalO(k),\Ker f )=0$ is
equivalent to $\Hom(\CalO(k),\CalF)=0$, which is an open
condition.
 The restricted maps $\pi_1$ and $\pi_2$ are still affine bundles of relative
 dimension $lr$ (see lemma \ref{buntor}).\\
We then have the correspondence 
\[
Irr(\underline{\Lambda}^{r,d,l}_{k,n,0}) \leftrightarrow
Irr(T^*\Vec_{k,0,0}^{r,d-l}) \times
Irr(\underline{\Lambda}^{0,l}_{k,n,0}).
\]
which preserves dimensions.\\
But as $T^*\Vec_{k,0,0}^{r,d-l}$ is an open substack of
$T^*\Vec^{r,d-l}$, which itself is an open substack of
$\underline{\Lambda}^{(r,d-l)}$, we have a natural inclusion of sets
\[
Irr(T^*\Vec_{k,0,0}^{r,d-l}) \subseteq
Irr(\underline{\Lambda}^{(r,d-l)}).
\]
 But we just proved that the
irreducible components on the right-hand side are of dimension $-r^2$. So
$T^*\Vec_{k,0,0}^{r,d-l}$ is pure of dimension $-r^2$, and by the
correspondence $\underline{\Lambda}^{r,d,l}_{k,n,0}$ is pure of
dimension $-r^2$. Suming up over the different values of $d,l$, we
obtain the desired result.
\begin{rem}
We may obtain the same result using section \ref{Irredcomp}. But this
proof may be generalized, for instance to weighted projective lines,
where a description of irreducible component is not known.
\end{rem}
\end{dem}
\subsection{Loop crystal operators}
We use theorem \ref{CorrespIrr} to define a combinatorial data similar
to a crystal graph on the set Irr($\underline{\Lambda}_\P1$).\\
Denote by $e_k^{max}$ the application
Irr($\underline{\Lambda}^{(r,d)}_{k,s,n}) \rightarrow
\text{Irr}(\underline{\Lambda}^{(r-s,d-sk)}_{k,0,n-s})$ (for $s>0$) deduced from Theorem \ref{CorrespIrr},
and $f_k^s$ its inverse.\\
Define the applications:
\[
f_k: \bigcup_{(r,d)}\text{Irr}(\underline{\Lambda}_\P1^{(r,d)}) \rightarrow
\bigcup_{(r,d)}\text{Irr}(\underline{\Lambda}_\P1^{(r-1,d-k)}) \cup \{0\}
\]
and
\[
e_k:\bigcup_{(r,d)}\text{Irr}(\underline{\Lambda}_\P1^{(r,d)}) \rightarrow
\bigcup_{(r,d)}\text{Irr}(\underline{\Lambda}_\P1^{(r+1,d+k)})
\]
as follows:\\
by Theorem \ref{dim2} if $Z\in\text{Irr}(\underline{\Lambda}_\P1^{(r,d)})$, there is a unique
$(s,n)$ such that $Z'=Z \cap \underline{\Lambda}^{(r,d)}_{k,s,n}$
is dense in $Z$. Then if $s>0$ we define
$f_k(Z)=\overline{e_k^{s-1}(f_k^{max}(Z'))}$, and $f_k(Z)=0$ if $s=0$. The map
$e_k$ is defined as $e_k(Z)=\overline{e_k^{s+1}(f_k^{max}(Z'))}$.\\
Let us define a map 
\[
\epsilon_k: \text{Irr}(\underline{\Lambda}) \rightarrow  \Nn
\]
by taking the generic value of the function $(\CalF,f) \mapsto -\text{rk}_k(\Ker f)$ on an irreducible component $Z$.\\
We have the map 
\[
wt: \text{Irr}(\underline{\Lambda}) \rightarrow \widehat{Q}
\]
corresponding to the decomposition $\text{Irr}(\underline{\Lambda})=\sqcup_\alpha \text{Irr}(\underline{\Lambda}^\alpha)$.\\
Now define the map 
\[
\begin{array}{cccc}
\phi_k: & \text{Irr}(\underline{\Lambda})  & \rightarrow & \Zz\\
 & Z &\mapsto & \epsilon_k(Z) + \langle h_1+k\delta,wt(Z) \rangle
\end{array}
\]
where $h_1$ is the fundamental weight in $P$.
\newline
\\
\textbf{Definition}: the set of data $\widehat{B}(\infty)=(\text{Irr}(\underline{\Lambda}_\P1),e_k,f_k,wt,\epsilon_k,\phi_k)$ is called the loop crystal associated to $\P1$.
\newline
\\
These data satisfy similar properties than crystals. For instance we
have the following:
\begin{enumerate}
\item $wt(f_k(Z))=wt(Z)-(\alpha_1+k\delta)$
\item $wt(e_k(Z))=wt(Z)+(\alpha_1+k\delta)$
\item $\epsilon_k(f_k(Z))=\epsilon_k(Z)+1$,
  $\epsilon_k(e_k(Z))=\epsilon_k(Z)-1$
\item $\phi_k(f_k(Z))=\phi_k(Z)-1$,
  $\phi_k(e_k(Z))=\phi_k(Z)+1$
\item If $f_k(Z)=Z'\neq 0$, then $e_k(Z')=Z$.
\end{enumerate}
The main difference is that we have now operators for any positive
root $\alpha_1+k\delta$.\\
As in the case of crystals, we can associate to this data a colored graph: vertices are indexed by elements of Irr($\underline{\Lambda}_\P1$) and $k$-colored edges are deduced from operators $f_k$.
\begin{prop}\label{connected}
The graph associated to the loop crystal $\widehat{B}(\infty)$ is connected.
\end{prop}
\begin{dem}
We divide the proof in two steps. The first step is the following
lemma:
\begin{lem}\label{connected1}
Any irreducible component $Z\in\text{Irr}(\underline{\Lambda})$ is
connected to some irreducible component $Z'$ with $\text{rk}(Z')=0$.
\end{lem}
\begin{dem}
We may consider that $\text{rk}(Z)>0$.\\
First we prove that $\lim_{k\rightarrow -\infty}\text{rk}_k(Z)>0$. Indeed
denote $Z=(\mathcal{V},\lambda)$, where $\mathcal{V}=\CalO(l_1)\oplus
\cdots \CalO(l_n)$ ($l_1\geq \cdots \geq l_n$) is a vector bundle
and $\lambda$ a partition.\\
On a generic point $(\CalF,f)$ in $Z$, we have $\CalF\simeq
\mathcal{V}\oplus \tau_\lambda$ and $\Im
(f)^\text{tor}=\tau_\lambda$. As $\Hom(\CalO(l_1),\CalO(k))=0$ for
$k<l_1$, we have $\Im(f_{|\CalO(l_1)})\subseteq \tau_\lambda$, then
$\text{rk}(\Ker f_{|\CalO(l_1)})=1$ and $\text{deg}(\Ker f_{|\CalO(l_1)})\geq
l_1-|\lambda|$. Then there exists an injection $\CalO(l_1-|\lambda|)
\hookrightarrow \Ker f$, which proves that $\text{rk}_k(Z)=\text{rk}_k(\CalF,f)>0$ for
$k\leq l_1-|\lambda|$.\\
Now for $Z$ take some $k<<0$ such that $\text{rk}_k(Z)>0$. Then
$f_k(Z)=Z'\neq 0$ with $\text{rk}(Z')=\text{rk}(Z)-1$. By easy
induction on $\text{rk}(Z)=n$, we can find
$k_1,\cdots k_n$ such that $f_{k_1}\circ \cdots \circ f_{k_n}(Z)$ is
an irreducible component of rank $0$.
\end{dem}
Now define $\mathbb{V}_k=\CalO(2k)\oplus \CalO(2(k-1)) \oplus \cdots \oplus
\CalO$. For a partition $\lambda=(\lambda_1,\cdots, \lambda_n)$ of
length $n$, set $\bar{\lambda}=(\lambda_1-1,\cdots ,\lambda_n-1)$,
a partition of length $l(\bar{\lambda})$ less or equal to $n$. We also fix $d=2k+2-l(\lambda)$. We prove the following:
\begin{lem}\label{connected2}
We have $f_d(\mathbb{V}_{k+1},\bar{\lambda})=(\mathbb{V}_k,\lambda)$.
\end{lem}
\begin{dem}
We work on generic parts of these irreducible components. A generic
point $(\CalF,f)\in Z=(\mathbb{V}_{k+1},\bar{\lambda})$ is such that
$\CalF \simeq \mathbb{V}_{k+1}\oplus \tau_{\bar{\lambda}}$ and $\text{rk}(\Ker f)$, $\text{deg}(\Ker f)$ are
minimal (or equivalently $\text{rk}(\Im f)$ and $\text{deg}(\Im f)$ are
maximal). The free part $(\Im f)^\text{fr}$ of the image of $f$ is then a subsheaf of
$\mathbb{V}_{k+1}(-2)=\mathbb{V}_k\oplus \CalO(-2)$ which is an image
of $\mathbb{V}_{k+1}$ of maximal rank. As there is no non zero
morphism from $\mathbb{V}_{k+1}$ to $\CalO(-2)$, $(\Im
f)^\text{fr}\subseteq \mathbb{V}_k$. But there is also no non zero
morphism from $\CalO(2k+2)$ to $\mathbb{V}_k$. Therefore $(\Im
f)^\text{fr}=\mathbb{V}_k$ and we may find a splitting $\mathbb{V}_k
\oplus \CalO(2k+2)$ such that 
\[
\begin{array}{cccc}\label{splitV}
f: & \CalO(2k+2) & \rightarrow & \tau_{\bar{\lambda}} \\
 & \mathbb{V}_k & \xrightarrow{\sim} & \mathbb{V}_k \\
& \tau_{\bar{\lambda}} & \rightarrow & \tau_{\bar{\lambda}}
\end{array}
\]
We claim that generically $(\Im f )^\text{tor} =\tau_{\bar{\lambda}}$. Indeed
, $\tau_{\bar{\lambda}}$ is
a sum of indecomposable sheaves with distinct support (generically on $Z$). Such a sheaf is isomorphic to the cokernel of a morphism $\CalO(j)\rightarrow
\CalO$ for some $j$. This proves that
there are surjections from any vector bundle to
$\tau_{\bar{\lambda}}$. So generically $(\Im f)^\text{tor}$ is isomorphic to $\tau_{\bar{\lambda}}$. \\
We proved that $\Im f =\mathbb{V}_k \oplus \tau_\lambda$. We then
deduce that $\text{rk}(\Ker f)=1$ and $\text{deg}(\Ker f)=2k+2$. But $\Ker
f^\text{tor}=\Ker f_{|\tau_{\bar{\lambda}}}$, and generically $\Ker
f^\text{tor}=soc(\tau_{\bar{\lambda}})$, the socle of the torsion
sheaf $\tau_{\bar{\lambda}}$ (which is of degree $l(\bar{\lambda})$).\\
Finally we have $\Ker f=soc(\tau_{\bar{\lambda}}) \oplus
\CalO(2k+2-l(\bar{\lambda}))=soc(\tau_{\bar{\lambda}}) \oplus
\CalO(d')$ if we set $d'=2k+2-l(\bar{\lambda})$. Then $\text{rk}_d(\Ker
f)=1$ since $d \leq d'$. Now take a generic injection $i:\CalO(d) \hookrightarrow \Ker
f$.\\
To see what the quotient $\CalF/i(\CalO(d))$ is, we use Serre's
description of coherent sheaves over $\P1$ as the category of graded
modules over the polynomial ring $\Cc[X,Y]$ modulo finite dimensional
modules. Define as usual $\Cc[X,Y]_j$
to be the graded $\Cc[X,Y]$-module with graduation shifted by $j$. If $x_j=(a_j:b_j)$ are the coordinates of the support of
$\tau_{\bar{\lambda}}$ in $\P1$, define $P(X,Y)=\Pi_j
(b_jX-a_jY)^{\bar{\lambda}_j}$ and $Q(X,Y)=\Pi_j(b_jX-a_jY)$. Now
the injection $i$ corresponds to a morphism:
\[
\begin{array}{cccc}
i:&\Cc[X,Y]_{d}& \rightarrow& \Cc[X,Y]_{d'}\oplus \Cc[X,Y]/Q\\
&1&\mapsto & (R(X,Y),\overline{S}(X,Y)),
\end{array}
\]
where $R$ is a homogeneous
polynomial of degree $d'-d$, $R$ and $S$ are prime to $P$ and $R$ has
$d'-d$ distinct roots (as we take a generic $i$). We will denote by
$y_j=(a_j':b_j')$ the roots of $R$. Then $y_j\neq x_k$ for any
$j,k$.\\
Because of \ref{splitV}, the injection $\Ker f \hookrightarrow \CalF$
factors through an injection $\Ker f \hookrightarrow \CalO(2k+2)
\oplus \tau_{\bar{\lambda}}$. In terms of $\Cc[X,Y]$-modules, this gives
rise to a morphism
\[
i: \Cc[X,Y]_{d'} \oplus \Cc[X,Y]/Q \rightarrow \Cc[X,Y]_{2k+2}\oplus
\Cc[X,Y]_{2k+2} \oplus \Cc[X,Y]/P
\]
given up to a constant by $i'(1,0)=(Q(X,Y),\overline{L}(X,Y))$ and
$i'(0,1)=(0,\overline{P/Q})$ for some $\overline{L}(X,Y)$. Generically, the
map $L:\CalO(d') \rightarrow \tau_{\bar{\lambda}}$ is surjective,
i.e. $L(X,Y)$ is prime to $P(X,Y)$.\\
Now the composition $i'\circ i$ is given by
\[
\begin{array}{cccc}
i'\circ i : & 1 & \mapsto & (R(X,Y)Q(X,Y),\overline{L}(X,Y)+\overline{P/Q}\overline{S}(X,Y))
\end{array}
\]
and we want to identify the cokernel of this morphism, which is the quotient $M$ of $\Cc[X,Y]_{2k+2}\oplus\Cc[X,Y]/P$ by the
submodule generated by the element:
\begin{equation}\label{generator}
(R(X,Y)Q(X,Y),\overline{P/Q}(X,Y)\overline{S}(X,Y)+\overline{L}(X,Y)).
\end{equation}
Now in this
quotient $M$ we have:
\begin{multline*}
Q(X,Y)R(X,Y)(1,0)=(R(X,Y)Q(X,Y),0)\\=-(0,\overline{L}(X,Y)+\overline{P/Q}\overline{S}(X,Y)
\end{multline*}
Observe that $(\overline{L}(X,Y)+\overline{P/Q}\overline{S}(X,Y))$ is prime to
$\overline{P}(X,Y)$, hence the submodule $N$ of $M$ generated by $(1,0)$ is of finite
codimension. Therefore $N\simeq M$ in $\Cc[X,Y]$-modgr/f.d.mod.\\
We claim that the coherent sheaf over $\P1$ associated to $N$ is a
torsion sheaf of
type $\lambda$. To see this, it suffices to determine the annihilator
$\text{Ann}_{(1,0)}$ of $(1,0)$. From \ref{generator} this annihilator is the
set of polynomials $A(X,Y)$ satisfying 
\begin{equation}\label{condideal}
(A(X,Y),0) \in \Cc[X,Y](RQ,\overline{P/Q} \overline{S} +\overline{L})
\end{equation}
In particular, we may write $A(X,Y)=RQB(X,Y)$ for some $B(X,Y) \in
\Cc[X,Y]$. Then \ref{condideal} becomes
\begin{equation*}
B(\overline{P/Q}\overline{S}+\overline{L})= 0 \text{  mod }P
\end{equation*}
But since $\overline{P/Q}\overline{S}+\overline{L}$ is prime to $P$, we see that
$B(X,Y) \in \Cc[X,Y]P(X,Y)$. We deduce that 
\[
\text{Ann}_{(1,0)}\simeq \Cc[X,Y]PQR
\]
and $N$ corresponds indeed to a sheaf of type $\lambda$.\\
Then $\CalF/i(\CalO(d))\simeq \mathbb{V}_k \oplus \tau_\lambda$ and
hence $f_d(Z)=(\mathbb{V}_k,\lambda)$ as wanted.
\end{dem}
We can now prove the proposition \ref{connected}: first use lemma
\ref{connected1} to reduce the problem to the connectedness between an
irreducible component $Z=(\emptyset,\lambda)$ of rank $0$ and
$\emptyset$. But by repeated use of lemma \ref{connected2}, we
easily see that $Z$ is connected to some $(\mathbb{V}_k,0)$. As
$f_{2k}\circ f_{2k-2} \circ \cdots \circ
f_0(\mathbb{V}_k,0)=\emptyset$, the result follows.
\end{dem}
This loop crystal is a kind of affine version of the crystal $B(\infty)$ for $sl_2$ contructed in \cite{Ka}. Moreover it carries a $\widehat{sl}_2$-crystal structure. To see define the following operators:
$\tilde{e}_1=e_0$, $\tilde{f}_1=f_0$, $\tilde{e}_0=f_{-1}$, $\tilde{f}_0=e_{-1}$, $\tilde{\epsilon}_1=\epsilon_0$, $\tilde{\phi}_1=\phi_0$, $\tilde{\epsilon}_0=\epsilon_{-1}$ and $\tilde{\phi}_0=\phi_{-1}$. Then we have:
\begin{prop}
The data $(\tilde{e}_i,\tilde{f}_i,wt,\tilde{\epsilon}_i,\tilde{\phi}_i, i=1,2)$ is an $\widehat{sl}_2$-crystal.
\end{prop}
We can see that this $\widehat{sl}_2$-structure is weaker than the original loop crystal structure as it is not connected.
\subsection{stability conditions}
This notion of loop crystal can be used to study representation of loop algebras. Following ideas from quiver varieties (see \cite{N}), crystals of interesting representations should be obtained by considering moduli spaces of \textit{stable} Higgs bundles with respect to an adequate notion of stability. For instance in the case of $\P1$ we can follow Laumon and Drinfeld (see \cite{La2}) and define the stable part of $\underline{\Lambda}_\P1$ as:
\[
\underline{\Lambda}_\P1^\text{s}=\{(\CalF,f)\in \underline{\Lambda}_\P1 \ | \ \Hom(\CalF,\CalF\otimes \Omega_\P1)^\text{nilp}=0\}
\]
When we consider this subcrystal as a $\widehat{sl}_2$-crystal and take the connected component of $\emptyset$ (= Irr($\underline{\Lambda}^0_\P1$)), the crystal obtained looks like a limit of crystals $B_n$, where $B_n$ is the (affine) crystal of the Kirillov-Reschetikhin module $V(n\varpi_1)$. We draw it below, denoting each irreducible component by the couple $(\mathcal{V},\lambda)$, where $\mathcal{V}$ is a vector bundle and $\lambda$ is a partition.
\begin{equation*}
\xymatrix{
\cdots \ar[r]^{f_0} & \CalO\oplus  \CalO \ar[r]^{f_0} & \CalO \ar[r]^{f_0} & \emptyset  \\
\cdots \ar[r]^{f_0} \ar[ru]^{f_{-1}} & \CalO\oplus  \CalO(-1)\ar[ru]^{f_{-1}} \ar[r]^{f_0} & \CalO(-1)  \ar[ru]^{f_{-1}} &  \\
\cdots \ar[r]^{f_0}\ar[ru]^{f_{-1}} & \CalO(-1)\oplus  \CalO(-1) \ar[ru]^{f_{-1}}&  &  \\
}
\end{equation*}
We can rewrite this by denoting an irreducible component by the rank $n$ and the degree $d$ (as there is only one such stable irreducible component of fixed rank and degree) by $n_d$. We only write the operators $\tilde{f}_1$ and $\tilde{f}_0$ of the $\widehat{sl}_2$-crystal.
\begin{equation*}
\xymatrix{
\cdots \ar @/^/ [r]^{\tilde{f}_1} & *+[F]{3_p} \ar @/^/ [r]^{\tilde{f}_1} \ar @/^/ [l]^{\tilde{f}_{0}[-1]} & 
 *+[F]{2_p} \ar @/^/ [r]^{\tilde{f}_1} \ar @/^/ [l]^{\tilde{f}_{0}[-1]}& *+[F]{1_p} \ar @/^/ [r]^{\tilde{f}_1} \ar @/^/ [l]^{\tilde{f}_{0}[-1]} & *+[F]{0_p} \ar @/^/ [l]^{\tilde{f}_{0}[-1]} 
}
\end{equation*}
The Kirillov-Reschetikhin modules are (conjecturally) the only finite dimensional modules which have a global basis (in the sense of Kashiwara). One can hope that more finite dimensional modules would have a "loop global basis" or a loop crystal. It would be interesting to find more stability conditions for arbitrary weighted projective lines and hence define global analogs of Nakajima's quiver varieties.


\begin{thebibliography}{LaMB}
\bibitem[BK]{BK} Bauman P., Kassel C., \textit{The Hall algebra of the
    category of coherent sheaves on the projective line}, J. reine
    angew. Math 533, 207-233 (2001).

\bibitem[Gr]{Gr} Grothendieck A., \textit{Techniques de construction
    et théorèmes d'existence en géométrie algébrique IV: les schémas
    de Hilbert}, séminaire Bourbaki, \textbf{221} (1960-61).

\bibitem[Jo]{Jo} Joyce D., \textit{Constructible functions on Artin
    stacks}, J. London. Math. Soc. (2) \textbf{74} No 3, 583--606 (2006).

\bibitem[Ka]{Ka} Kashiwara, M., \textit{On crystal bases of the
    $Q$-analogue of universal enveloping algebras}, Duke Math. J.,
  \textbf{63},  no. 2, 465--516  (1991). 

\bibitem[Kap]{Kap} Kapranov M., \textit{Eisenstein Series and quantum
    affine algebras}, J. Math. Sci. (New York) \textbf{84}, 1311--1360 (1997).

\bibitem[KS]{KS} Kashiwara M., Saito Y., \textit{Geometric
    construction of crystal bases}, Duke Math. J., \textbf{89}, 9-36 (1997).

\bibitem[La1]{La1} Laumon G., \textit{Correspondance de Langlands
    géométrique pour les corps de fonctions}, Duke
    Math. J., \textbf{54}, 309-359 (1987).

\bibitem[La2]{La2} Laumon G., \textit{Un analogue global du cône
    nilpotent}, Duke Math. J., \textbf{57}, 647-671 (1988).

\bibitem[LaMB]{LaMB} Laumon G., Moret-Bailly L., \textit{Champs
    Algébriques}, Ergeb. Math. Grenzgeb. 3.Folge, vol.39, Springer 2000.

\bibitem[Le]{Le} LePotier J., \textit{Fibrés vectoriels sur les
    courbes algébriques}, Publ.Math.Univ.Paris VII (1996).

\bibitem[L1]{L1} Lusztig G., \textit{Semicanonical bases arising from
    enveloping algebras}, Advances in Mathematics, \textbf{151},
    129-139 (2000).

\bibitem[L2]{L2} Lusztig G. \textit{Canonical bases arising from quantized enveloping algebras}, J. Amer. Math. Soc. \textbf{3}, 447-498 (1990).

\bibitem[N]{N} H. Nakajima, \emph{Quiver varieties and finite-dimensional representations of quantum affine algebras},  J. Amer. Math. Soc.  \textbf{14}  (2001),  no. 1, 145--238 

\bibitem[Ri]{Ri} Ringel C.M., \textit{Hall algebras and quantum groups}, Invent. Math., \textbf{101}, 583-591 (1990).

\bibitem[Sc1]{Sc1} Schiffmann O., \textit{Noncommutative projective curves and quantum loop algebras}, Duke Mat. J. \textbf{121}, 113-168 (2004).
\bibitem[Sc2]{Sc2} Schiffmann O., \textit{Canonical bases and moduli
    space of sheaves on curves}, Inventiones Mathematicae,
    \textbf{165}, 453-524 (2006).

\end{thebibliography}
\end{document}